\documentclass{article}
\title{Mukai flops and derived categories II}
\author{Yoshinori Namikawa}
\date{ }
\begin{document}
\maketitle

\begin{center}
{\bf Introduction}
\end{center}

This paper is a sequel to [Na]. 
Let $X$ and $Y$ be birationally 
equivalent smooth quasi-projective 
varieties. Then we say that $X$ and 
$Y$ are $K$-equivalent if there is a 
smooth quasi-projective variety $Z$ 
with proper birational morphisms   
$f : Z \to X$ and 
$g : Z \to Y$ such that $f^*K_X = 
g^*K_Y$.  The following problem is 
a motivation of this paper. 
\vspace{0.12cm}

{\bf Problem 1}. {\em Let 
$X$ and $Y$ be $K$-equivalent 
smooth quasi-projective varieties. 
Then, is there an equivalence of 
bounded derived 
categories of coherent sheaves   
$D(X) \to D(Y)$? }  
\vspace{0.12cm}

In some good cases, birational maps  
$X - - \to Y$ can be decomposed 
into certain kinds of flops.  A flop is 
a diagram of 
quasi-projective varieties 
$$ X \stackrel{s}\rightarrow {\bar X} 
\stackrel{s^+}\leftarrow X^+.$$
Here $s$ and $s^+$ are small, projective, 
crepant resolutions 
of the normal variety ${\bar X}$ with 
$\rho (X/{\bar X}) = \rho (X^+/{\bar X}) = 1$. 
Moreover, for an $s$-negative divisor $D$, its 
proper transform $D'$ becomes $s^+$-ample. 
A special case of Problem 1 is:   
\vspace{0.12cm}

{\bf Problem 2}. {\em  Let 
$$ X \stackrel{s}\rightarrow {\bar X} 
\stackrel{s^+}\leftarrow X^+$$ be a flop. 
Then, is there an equivalence of bounded derived 
categories of coherent sheaves 
$D(X) \to D(X^+)$? }   
\vspace{0.12cm}

If Problem 2 is affirmative, then, 
which functor gives the equivalence ? 
The following examples suggest that 
the functor $\Psi$ defined by the fiber product 
$X \times_{{\bar X}}X^+$ would be a 
correct one.  \vspace{0.12cm}

{\bf Examples}. (1) ([B-O]): Let $X$ be a 
smooth quasi-projective variety of dimension 
$2h-1$ which contains a subvariety $M 
\cong {\mathbf P}^{h-1}$ with $N_{M/X} 
\cong {\mathcal O}(-1)^{\oplus h}$. One can blow up 
$X$ along $M$ and blow down the exceptional 
divisor in another direction. In this way, we have a 
new variety $X^+$ with a subvariety $M^+ 
\cong {\mathbf P}^{h-1}$. Let $s: X \to {\bar X}$ 
and $s^+: X^+ \to {\bar X}$ be the birational 
contraction maps of $M$ and $M^+$ to points 
respectively.  Let $\mu: X \times_{{\bar X}}X^+ 
\to X$ and $\mu^+: X \times_{{\bar X}}X^+ \to 
X^+$ be the projections. Then $\Psi (\bullet) 
:= {\mathbf R}{\mu^+}_*{\mathbf L}\mu^*(\bullet)$ 
is an equivalence.     
\vspace{0.12cm}

(2) ([Na], [Ka 2]): Let $X$ be a smooth 
quasi-projective variety of dimension 
$2h-2$ which contains a subvariety $M 
\cong {\mathbf P}^{h-1}$ with $N_{M/X} 
\cong \Omega^1_{{\mathbf P}^{h-1}}$. Let 
$s: X \to {\bar X}$ be a birational contraction 
map of $M$ to a point. We have a flop 
$$ X \stackrel{s}\rightarrow {\bar X} 
\stackrel{s^+}\leftarrow X^+$$ 
along $M$. Then the functor $\Psi$ defined 
by the fiber product $X \times_{{\bar X}}X^+$ 
gives an equivalence $D(X) \to D(X^+)$.    
\vspace{0.12cm}

(3) ([Br 1]; cf. also [Ch],[Ka 2],[vB]):  Let 
$$ X \stackrel{s}\rightarrow 
{\bar X} \stackrel{s^+}\leftarrow 
X^+$$ be a 3-dimensional flop with 
$X$ and $X^+$ being smooth quasi-projective 
3-folds with trivial canonical line bundles.  
Then the functor $\Psi$ defined 
by the fiber product $X \times_{{\bar X}}X^+$ 
gives an equivalence $D(X) \to D(X^+)$.   
\vspace{0.15cm}

Markman has studied in [Ma] a generalization 
of the Mukai flop. Here we call it a 
{\em stratified Mukai flop}. 
In this paper we observe that, for a 
stratified Mukai flop, \vspace{0.12cm}     

(1) the fiber product $X \times_{{\bar X}}X^+$ 
defines an isomorphism $$K(X) \cong K(X^+)$$ 
of Grothendieck groups, but 

(2) the functor $\Psi$ is not an equivalence in 
general.   
\vspace{0.12cm}
   
More precisely, let $H$ be a ${\mathbf C}$-vector 
space of $\dim h$. For $t \leq h/2$, let 
$G:= G(t,H)$ be the Grassmann variety parametrizing 
$t$ dimensional vector subspaces of $H$. 
We denote by $T^*G$ the cotangent bundle of $G$. 
The nilpotent variety ${\bar N}^t(H)$ is defined as  
$${\bar N}^t(H) := \{A \in \mathrm{End}(H); A^2 = 0, 
\mathrm{rank}(A) \leq t\}.$$    
Then there is a natural morphism $s: T^*G \to 
{\bar N}^t(H)$ and this gives a resolution of the 
singular variety ${\bar N}^t(H)$. Let $H^*$ be the 
dual space of $H$ and put $G^+ := G(t, H^*)$. 
Now we have the resolution $s^+: T^*G^+ \to 
{\bar N}^t(H^*)$. Since there is a 
natural identification ${\bar N}^t(H) \cong 
{\bar N}^t(H^*)$, we have a commutative 
diagram 
$$ T^*G - - \to T^*G^+ $$ 
$$ \downarrow \hspace{1.2cm} \downarrow $$ 
$$ {\bar N}^t(H) \cong {\bar N}^t(H^*).$$ 

Moreover, we define 1-parameter deformations 
$E(H) \to {\mathbf C}^1$ of $T^*G$ and 
${\bar{\mathcal N}}^t(H) \to {\mathbf C}^1$ of 
${\bar N}^t(H)$.  The map $s: T^*G \to 
{\bar N}^t(H)$ extends to ${\tilde s}:  E(H) 
\to {\bar{\mathcal N}}^t(H)$. For the dual space 
$H^*$ we construct the same, and 
we have a commutative diagram \newpage
$$ E(H) - - \to E(H^*) $$ 
$$ \downarrow \hspace{1.2cm} \downarrow $$ 
$$ {\bar{\mathcal N}}^t(H) \cong {\bar{\mathcal N}}^t(H^*).$$ 
When $t = 1$, these diagrams give Example (2) 
and Example (1) respectively. We call the first diagram 
a {\em stratified Mukai flop}\footnote{This notation is not exact. 
In fact, when $t = h/2$, 
the birational map $T^*G --\to T^*G^+$ becomes an isomorphism. 
But we include here such cases because our functor $\Psi_0$ 
is not a trivial one induced by the isomorphism} 
 and call the second diagram 
a {\em stratified Atiyah flop}.  

Let $\Psi_0: D(T^*G^+) \to 
D(T^*G)$ be the functor defined by 
$T^*G \times_{{\bar N}^t(H)}T^*G^+ \in D(T^*G \times 
T^*G^+)$.  
Similarly, let 
$\Psi : D(E(H^*)) \to D(E(H))$ be the 
functor defined by 
$E(H) \times_{{\bar{\mathcal N}}^t(H)} E(H^*) \in 
D(E(H)\times E(H^*))$.  
These functor naturally induces homomorphisms 
of Grothendieck groups. Therefore, we have two 
commutative diagrams 
$$  D(T^*G^+) \stackrel{\Psi_0}\to D(T^*G) $$
$$  \downarrow \hspace{1.2cm} \downarrow $$ 
$$  K(T^*G^+) \stackrel{{\bar \Psi}_0}\to 
K(T^*G).$$   

$$  D(E(H^*)) \stackrel{\Psi}\to D(E(H)) $$
$$  \downarrow \hspace{1.2cm} \downarrow $$ 
$$  K(E(H^*)) \stackrel{{\bar \Psi}}\to K(E(H)).$$  

We shall prove that both ${\bar \Psi}$ 
and ${\bar \Psi}_0$ are isomorphisms (Theorem (2.6)
and Theorem (2.7)). 
As noted above, when $t = 1$, 
$\Psi_0$ and $\Psi$ are both equivalences. 
But, when $t = 2$ and $h = 4$, none of $\Psi_0$ 
or $\Psi$ is equivalence (Section 4). To show this, 
we have to describe the birational map 
$E(H) - - \to E(H^*)$ as the composite of  
certain blowing-ups and blowing-downs [Ma]. 
In section 3, we shall sketch this when $t = 2$ 
and $h = 4$.  
At the moment, when $t \geq 2$, the following 
question is completely open: 
\vspace{0.12cm}

{\bf Problem 3}. {\em For a generalized Mukai flop 
(resp. a generalized Atiyah flop), does the isomorphism 
${\bar \Psi}_0$ (resp. ${\bar \Psi}$) lift to an 
equivalence of derived categories ?}
\vspace{0.12cm}  
   
The above constructions can be also applied to 
flag varieties. The final section deals with the case 
of complete flag varieties (cf. [Sl]). 
In this case, we have an 
equivalence of the bounded derived categories of 
coherent sheaves on dual pairs (Theorem (5.9.1)).        
\vspace{0.15cm}\newpage
    
\begin{center}
{\bf \S 1.} {\bf Stratified Mukai flops}
\end{center}

(1.1). Let $H$ be a ${\mathbf C}$-vector space of 
dimension $h$. For a positive integer $t \leq 
h/2$, let $G(t, H)$ be the Grassmann variety 
parametrizing $t$ dimensional vector subspaces 
of $H$. In the remainder we simply write $G$ for 
$G(t, H)$. Let $\tau$ and $q$ be the universal 
subbundle and the universal quotient bundle 
respectively. They fits into the exact sequence 
$$ 0 \to \tau \to H\otimes_{\mathbf C}{\mathcal O}_G 
\to q \to 0.$$

(1.2). Let $T^*G$ be the cotangent bundle of $G$ 
and let $\pi: T^*G \to G$ be the projection map. 
The nilpotent variety ${\bar N}^t(H)$ is defined 
as 
$$ {\bar N}^t(H) := \{A \in \mathrm{End}(H); 
A^2 = 0, \mathrm{rank}(A) \leq t\}. $$ 

We define a birational morphism 
$$s: T^*G \to {\bar N}^t(H)$$ 
as follows.  
First note that $T^*G$ is identified with 
the vector bundle ${\underline{\mathrm{Hom}}}(q, \tau)$ 
over $G$. Then a point of $T^*G$ is expressed 
as a pair $(p, \phi)$ of a point $p \in G$ and 
an element $\phi \in \mathrm{Hom}(q(p), \tau(p))$. 
Since there is a surjection $H \to q(p)$ and an 
injection $\tau(p) \to H$, the element $\phi$ 
defines an element of $\mathrm{End}(H)$. We 
denote this element by the same $\phi$. 
Now we define $s((p, \phi)) := \phi$. 
It is easy to check that $\phi \in  
{\bar N}^t(H).$  We see that, for $A \in 
{\bar N}^t(H)$, 
$$s^{-1}(A) = \{(p, A); \mathrm{Im}(A) 
\subset \tau(p) \subset \mathrm{Ker}(A)\}.$$  
For $i \leq t$, put   
$N^i(H) := \{A \in \mathrm{End}(H); A^2 = 0, 
\mathrm{rank}(A) = i\}$. Then, for $A \in 
N^i(H)$,  
$$s^{-1}(A) \cong G(t-i,h-2i).$$ 
In particular, $s$ is an isomorphism 
over $N^t(H)$.  
The map $s$ is called the 
Springer resolution of ${\bar N}^t(H)$.

(1.3). We construct a 1-parameter deformation 
of the Springer resolution $s: T^*G \to 
{\bar N}^t(H)$. First of all, we shall define 
a vector bundle $E(H)$ over $G$ and an exact sequence 
$$ 0 \to {\underline{\mathrm{Hom}}}(q, \tau) \to E(H) 
\to {\mathcal O}_G \to 0. $$ 
As in (1.2), $\mathrm{Hom}(q(p), \tau(p))$ 
is embedded in $\mathrm{End}(H)$. For 
a suitable basis of $H$,   
$\mathrm{Hom}(q(p), \tau(p))$ is the 
set of $h \times h$ matrices of the   
following form: 

$$\left(
\begin{array}{cc}
0 & *\\
0 & 0  
\end{array}\right),$$ 
where the first $0$ is the $t \times t$ 
zero matrix and where * is a $t \times 
(h-t)$ matrix.  
Now, let $E(H)(p)$ be the set of 
$h \times h$ matrices which have the 
following form: 
    
$$\left(
\begin{array}{cc}
\alpha I & *\\
0 & 0  
\end{array}\right),$$ 
where $\alpha \in {\mathbf C}$ 
and where $I$ is the $t \times t$ 
identity matrix. $E(H)(p)$ is characterized 
as the set of elements $\phi \in {\mathrm 
End}(H)$ such that $\mathrm{Im}(\phi) 
\subset \tau(p)$ and $\phi\vert_{\tau(p)} 
= \alpha I$; hence it is independent of 
the choice of the basis of $H$.   
By the construction, there is an 
exact sequence 
$$ 0 \to \mathrm{Hom}(q(p), \tau(p)) 
\to E(H)(p) \to {\mathbf C} \to 0. $$
Note that the map 
$E(H)(p) \to {\mathbf C}$ 
is defined as $(1/t){\mathrm trace}$; 
hence it is also independent of the 
choice of the basis of $H$.  
We put $E(H) := \cup_{p \in G}E(H)(p)$. 
Then $E(H)$ becomes a vector bundle over 
$G$, and there is an exact sequence of 
vector bundles 
$$ 0 \to {\underline{\mathrm{Hom}}}(q, \tau) \to E(H) 
\to {\mathcal O}_G \to 0. $$ 
There is a surjective morphism from $E(H)$ to 
${\mathbf C}^1$, and its central fiber is 
$T^*G$.   
 
(1.4). Define  
$${\mathcal N}^t(H) := \{A \in \mathrm{End}(H); 
A^2 = (\exists{\mathrm scalar})A , 
\mathrm{rank}(A) = t \}. $$ 
Every elements of 
${\mathcal N}^t(H)$ is conjugate to a matrix 
 
$$\left(
\begin{array}{cc}
\alpha I & *\\
0 & 0  
\end{array}\right)$$ 
where $\alpha I$ is a $t \times t$ 
scalar matrix and where * is a $t \times 
(h-t)$ matrix. 
Now let ${\bar{\mathcal N}}^t(H)$ be the 
Zariski closure of ${\mathcal N}^t(H)$ in 
$\mathrm{End}(H)$. By taking $(1/t){\mathrm 
trace}$, we define a morphism  
${\bar{\mathcal N}}^t(H) \to {\mathbf C}^1$. 
Its central fiber is ${\bar N}^t(H)$.  
  
(1.5). Each point of $E(H)$ is expressed as 
a pair $(p, \phi)$ of $p \in G$ and $\phi 
\in E(H)(p)$. The extended Springer 
resolution 
$$ {\tilde s}: E(H) \to {\bar{\mathcal N}}^t(H)$$ 
is defined as ${\tilde s}((p, \phi)) := \phi$.  

As a consequence, we have a 1-parameter 
deformation 
$$ E(H) \to {\bar{\mathcal N}}^t(H) \to {\mathbf C}^1$$ 
of 
$$ T^*G \to {\bar N}^t(H) \to 0. $$  
 
(1.6). We shall define {\em dual objects} 
of those constructed above. Let $H^*$ be the 
dual space of $H$. For $H^*$, we define the similar 
objects to (1.1), ..., (1.4). We write $G^+$ for 
$G(t, H^*)$. $s^+: T^*G^+ \to {\bar N}^t(H)$ and 
${\tilde s}^+ : E(H^*) \to {\bar{\mathcal N}}^t(H^*)$ are 
the Springer resolution and the extended Springer 
resolution, respectively. The relationship between 
these duals are as follows. There is a canonical 
isomorphism $\mathrm{End}(H) \cong 
\mathrm{End}(H^*)$. With respect to dual bases of 
$H$ and $H^*$, this isomorphism is given by 
the transpose $A \to {}^tA$. This isomorphism 
naturally identifies ${\bar N}^t(H)$ with 
${\bar N}^t(H^*)$, and ${\bar{\mathcal N}}^t(H)$ 
with ${\bar{\mathcal N}}^t(H^*)$, respectively. 
These identifications induce a birational map 
between $T^*G$ and $T^*G^+$, and a birational 
map between 
$E(H)$ and $E(H^*)$:      

$$ T^*G  - - \to T^*G^+ $$ 
$$ \downarrow \hspace{1.5cm} \downarrow $$ 
$$ {\bar N}^t(H) \hspace{0.3cm} \cong 
\hspace{0.3cm} {\bar N}^t(H^*),$$ 

$$ E(H) - - \to E(H^*) $$ 
$$ \downarrow \hspace{1.5cm} \downarrow $$ 
$$ {\bar{\mathcal N}}^t(H) \hspace{0.3cm} 
\cong \hspace{0.3cm} {\bar{\mathcal N}}^t(H^*).$$ 

{\bf Example (1.7.1)}: When $t = 1$ and 
$h = 2$, the birational map $T^*G - - \to 
T^*G^+$ is an isomorphism. 
$T^*G$ (resp. $T^*G^+$) is a non-singular 
surface and its zero section is a $(-2)$-curve. 
The nilpotent 
variety ${\bar N}:= {\bar N}^t(H) \cong 
{\bar N}^t(H^*)$ has an $A_1$ surface 
singularity at $0 \in {\bar N}$. The Springer 
resolution $s$ (or $s^+$) is nothing but the 
minimal resolution of ${\bar N}$.  
On the other hand, $E(H)$ and $E(H^*)$ are 
both non-singular 3-folds, and their zero sections  
are $(-1,-1)$-rational curves. 
The birational map $E(H) - - \to E(H^*)$ 
is the {\em Atiyah flop} along these 
$(-1,-1)$-curves. ${\bar{\mathcal N}}:= 
{\bar{\mathcal N}}^t(H) \cong {\bar{\mathcal N}}^t
(H^*)$ has an ordinary double point at 
the origin. The extended Springer resolutions 
${\tilde s}$ and ${\tilde s}^+$ are mutually 
different small resolutions of this ordinary 
double point.   
\vspace{0.12cm}

{\bf Example (1.7.2)}: When $t = 1$ and 
$h \geq 3$, $T^*G$ and $T^*G^+$ are non-singular 
varieties of dim $2(h-1)$. The zero sections 
of them are isomorphic to ${\mathbf P}^{h-1}$. 
The birational map $T^*G - - \to T^*G^+$ 
is the {\em Mukai flop} along these ${\mathbf 
P}^{h-1}$. On the other hand, $E(H)$ and 
$E(H^*)$ are non-singular varieties of dim 
$2h -1$. The zero sections of them are 
${\mathbf P}^{h-1}$ whose normal bundles are 
isomorphic to ${\mathcal O}(-1)^{\oplus h}$.  
The birational map $E(H) - - \to E(H^*)$ is 
the flop treated in [B-O].     
\vspace{0.15cm}

\begin{center}
{\bf \S 2.} {\bf K-theory}
\end{center}

(2.1). For an algebraic scheme $X$, we denote 
by $K(X)$ the Grothendieck group of $X$. Namely, 
let $F(X)$ be the free Abelian group 
whose basis consists of the set of coherent 
sheaves on $X$. Let 
$$ (E): 0 \to {\mathcal F}' \to {\mathcal F} 
\to {\mathcal F}'' \to 0 $$
be an exact sequence of coherent sheaves on 
$X$. The exact sequence $(E)$ defines an 
element $Q(E) := [{\mathcal F}] - [{\mathcal F}'] 
- [{\mathcal F}'']$ of $F(X)$. $K(X)$ is the 
quotient of $F(X)$ by the subgroup 
generated by all such $Q(E)$. 

(2.2).  Let $G:= G(t, H)$ be the same as (1.1).
Let $\alpha := (\alpha_1, ..., \alpha_n)$ be a 
sequence of non-negative integers with 
$\alpha_1 \geq \alpha_2, ...,  \geq \alpha_n$. 
One can associate a Young diagram with $\alpha$. 
Denote by $r(\alpha)$ the number of the rows 
of this Young diagram, and denote by $c(\alpha)$ 
the number of the columns of this Young diagram. 
For such an $\alpha$ and for a vector bundle 
$E$ of rank $n$ over $G$, one can define a new vector bundle 
$\Sigma^{\alpha}E$ over $G$ (cf.[Kap],[Fu 1, {\S}.8.]).  
The following 
theorem is well-known.  \vspace{0.15cm}

{\bf Theorem (2.2.1)}. {\em $K(G)$ is a free abelian 
group generated by $[\Sigma^{\alpha}\tau]$ 
with  $r(\alpha) \leq t$ and $c(\alpha) \leq h-t$.}  
\vspace{0.12cm}

{\em Sketch of the Proof}. 
By [Kap], the bounded derived category 
$D(G)$ of coherent sheaves is generated by 
$[\Sigma^{\alpha}\tau]$ with 
$r(\alpha) \leq t$ and $c(\alpha) \leq h-t$. 
In particular, $K(G)$ is 
generated by them. On the other hand, the Chern 
character homomorphism $K(G) \to 
A(G)_{{\mathbf Q}}$ induces an isomorphism 
$K(G)_{{\mathbf Q}} \cong A(G)_{{\mathbf Q}}$ 
(cf. [Fu 2, Example 15.2.16]).  
Moreover, $A(G)$ is a free ${\mathbf Z}$-module 
with the basis $\{\alpha\}$, the Schubert classes 
for Young tableaux $\alpha$ with 
$r(\alpha) \leq t$ and $c(\alpha) \leq h-t$ 
(cf.[Fu 2, 14.7]). 
\vspace{0.15cm}
 
{\bf Example (2.2.2)}. $K({\mathbf P}^{h-1})$ is 
generated by ${\mathcal O}$, ${\mathcal O}(-1)$, ...., 
and ${\mathcal O}(-h+1)$.                                            
\vspace{0.12cm}

Let $\pi: T^*G \to G$ and ${\tilde \pi}: E(H) 
\to G$ be the same as Section 1.  Since 
they are vector bundles over $G$, the  
natural maps $\pi^* : K(G) \to K(T^*G)$ 
and ${\tilde \pi}^*: K(G) \to K(E(H))$ are 
both isomorphisms (cf. [C-G, Theorem 5.4.17]). 
Write $\tau_{T^*G}$ for ${\pi}^*\tau$ and 
$\tau_{E(H)}$ for ${\tilde \pi}^*\tau$. 
\vspace{0.15cm}

{\bf Corollary (2.2.3)}. (1) {\em $K(T^*G)$ is 
a free Abelian group generated by 
$[\Sigma^{\alpha}\tau_{T^*G}]$ with $r(\alpha) 
\leq t$ and $c(\alpha) \leq h-t$.} 

(2) {\em $K(E(H))$ is a free Abelian group 
generated by $[\Sigma^{\alpha}\tau_{E(H)}]$ 
with $r(\alpha) \leq t$ and $c(\alpha) \leq h-t$.}     
\vspace{0.15cm}

(2.3). Let $\psi: E(H) - - \to E(H^*)$ be the birational 
map in (1.6). Note that the extended Springer resolution 
${\tilde s}: E(H) \to {\bar{\mathcal N}}^t(H)$ (resp. 
${\tilde s}^+: E(H^*) \to {\bar{\mathcal N}}^t(H^*)$) is 
a small resolution; hence $\psi$ is 
an isomorphism in codimension 1. Namely, there are 
Zariski open subsets $U \subset E(H)$ and $U^+ \subset 
E(H^*)$ such that $U \cong U^+$ and the complement of 
$U$ in $E(H)$(resp. $U^+$ in $E(H^*)$) is at least of 
codimension 2. Let $F$ be a reflexive sheaf on $E(H)$. 
Since $U \cong U^+$, $F\vert_U$ is regarded as a sheaf 
on $U^+$. Then the direct image $\psi(F) := (j^+)_*(F\vert_U)$ 
under the inclusion map $j^+: U^+ \to E(H^*)$ is a 
reflexive sheaf on $E(H^+)$. We call $\psi(F)$ the proper 
transform of $F$.  \vspace{0.15cm}

{\bf Lemma (2.3.1)}. {\em $\psi(\tau_{E(H)}) = 
({\tau^+}_{E(H^*)})^*$, where 
$({\tau^+}_{E(H^*)})^*$ is the dual sheaf of ${\tau^+}_{E(H^*)}$.}       
\vspace{0.12cm}

{\em Proof}. Put $M := 
{\bar{\mathcal N}}^t(H) - \mathrm{Sing}({\bar{\mathcal N}}^t(H))$ 
and  
$M^+ := {\bar{\mathcal N}}^t(H^*) 
- \mathrm{Sing}({\bar{\mathcal N}}^t(H^*))$. 
As explained in (1.6), there is 
an isomorphism, say $\psi_0 : M 
\cong M^+$. For $A \in M$, $(\tau_{E(H)})(A) = 
\mathrm{Im}(A) 
\subset H$. Since $\psi_0(A) = { }^tA$ (cf. (1.6)), 
$((\psi_0)_*{\tau}_{E(H)})({}^tA) = \mathrm{Im}(A)$. In other 
 words, for $B \in M^+$, $((\psi_0)_*{\tau}_{E(H)})(B) = 
\mathrm{Im}({}^tB)$. Note that $\tau^+_{E(H^*)}(B) = {\mathrm 
im}(B)$, where $B$ is an endomorphism of $H^*$. 
Therefore, $({\tau^+}_{E(H^*)})(B)^* = \mathrm{Im}({}^tB)$. 
This implies that             
$((\psi_0)_*\tau_{E(H)})(B) = ({\tau^+}_{E(H^*)})(B)^*$.  
Since $M$ and $M^+$ are naturally embedded as 
open subsets of $E(H)$ and $E(H^*)$ and their 
complements are of codimension at least 2, we 
conclude that $\psi(\tau_{E(H)}) = ({\tau^+}_{E(H^*)})^*$.        
\vspace{0.15cm}

{\bf Corollary (2.3.2)}. {\em 
For an $\alpha$ of (2.2), 
 $\psi(\Sigma^{\alpha}
\tau_{E(H)}) = 
(\Sigma^{\alpha}
\tau^+_{E(H^*)})^*.$}   
\vspace{0.15cm}

(2.4). Recall that ${\bar N} := 
{\bar N}^t(H) \cong {\bar N}^t(H^*)$ 
and ${\bar{\mathcal N}} := 
{\bar{\mathcal N}}^t(H) \cong 
{\bar{\mathcal N}}^t(H^*)$. 
Let $\mu_0: T^*G \times_{{\bar N}} 
T^*G^+ \to T^*G$ and 
$\mu^+_0: T^*G \times_{{\bar N}} 
T^*G^+ \to T^*G^+$ be the natural 
projections. Similarly, let 
$\mu: E(H) \times_{{\bar{\mathcal N}}} 
E(H^*) \to E(H)$ and 
$\mu^+: E(H) \times_{{\bar{\mathcal N}}} 
E(H^*) \to E(H^+)$ be the natural 
projections.     
We define a homomorphism 
$$\Psi_0: K(T^*G^+) \to K(T^*G)$$ 
as $\Psi_0(\bullet) :=  (\mu_0)_*\circ 
(\mu^+_0)^*(\bullet)$.
Similarly, a homomorphism  
$$\Psi: K(E(H^*)) \to K(E(H))$$ 
is defined as $\Psi(\bullet) := 
\mu_*\circ (\mu^+)^*(\bullet)$. 
\vspace{0.12cm}

{\bf Proposition (2.5)}. {\em For 
generators $[\Sigma^{\alpha}
\tau^+_{E(H^*)}]$ of $K(E(H^*))$
(cf. Corollary (2.2.3)), 
$\Psi([\Sigma^{\alpha}
\tau^+_{E(H^*)}]) = 
[\Sigma^{\alpha}
(\tau_{E(H)})^*].$}   
\vspace{0.12cm}

{\em Proof}. Let $F^+$ be one of 
$\Sigma^{\alpha}
\tau^+_{E(H^*)}$'s.  
By definition, 
$$\Psi([F^+]) 
= \Sigma (-1)^i [R^i\mu_*({\mu^+}^*F^+)].$$ 
Note that, for $i > 0$, 
$R^i\mu_*({\mu^+}^*F^+)$ has the 
support in $T^*G$. The cokernel 
of the injection: 
$$\mu_*({\mu^+}^*F^+) 
\to  (\mu_*({\mu^+}^*F^+))^{**}$$ 
also has the support in $T^*G$. 
Here, $**$ 
means the double dual.
Hence, 
by the following lemma, $[R^i\mu_*({\mu^+}^*F^+)] 
= 0$ if $i > 0$, and $[\mu_*({\mu^+}^*F^+)] 
= [(\mu_*({\mu^+}^*F^+))^{**}]$.  By Corollary (2.3.2), 
$(\mu_*({\mu^+}^*F^+) )^{**} \cong 
\Sigma^{\alpha}(\tau_{E(H)})^*$.  
\vspace{0.12cm}

{\bf Lemma (2.5.1)}. {\em Let $F$ be a coherent 
sheaf on $E(H)$ whose support is contained in 
$T^*G$. Then $[F] = 0$ in $K(E(H))$.}    
\vspace{0.12cm}

{\em Proof}. Let $I$ be the ideal sheaf of 
$T^*G$ in $E(H)$. For a sufficiently large 
$n > 0$, $I^nF = 0$. Hence, $F$ can be described 
as successive extensions of ${\mathcal O}_{T^*G}$-modules.  
It suffices to show that  every 
coherent ${\mathcal O}_{T^*G}$ module is 
the zero as an element of $K(E(H))$.  
Since $K(T^*G)$ is generated by 
the elements of the form:  
$[\Sigma^{\alpha}
\tau_{T^*G}]$, we have to 
prove that  $[\Sigma^{\alpha}
\tau_{T^*G}] = 0$ in $K(E(H))$. 
Note here that $\Sigma^{\alpha}
\tau_{T^*G} = \Sigma^{\alpha}
\tau_{E(H)}\vert_{T^*G}$.   
Since $E(H)$ is a 1-parameter 
deformation of $T^*G$, there is an 
exact sequence 
$$ 0 \to {\mathcal O}_{E(H)} \to 
{\mathcal O}_{E(H)} \to {\mathcal O}_{T^*G} 
\to 0.$$
By taking the tensor product 
with $\Sigma^{\alpha}
\tau_{E(H)}$, we have the exact 
sequence:  
$$ 0 \to \Sigma^{\alpha}
\tau_{E(H)} \to \Sigma^{\alpha}
\tau_{E(H)} \to \Sigma^{\alpha}
\tau_{T^*G} \to 0. $$    
This yields that $[\Sigma^{\alpha}
\tau_{T^*G}] = 0$ in $K(E(H))$.      
\vspace{0.15cm}

{\bf Theorem (2.6)}. {\em 
The homomorphism  
$$\Psi: K(E(H^*)) \to K(E(H))$$ 
is an isomorphism.}  
\vspace{0.12cm}

{\em Proof}. By (2.5), it is sufficient to prove that 
$[\Sigma^{\alpha}(\tau_{E(H)})^*]$ with 
$r(\alpha) \leq t$ and $c(\alpha) \leq h-t$ 
form a basis of $K(E(H))$. 
We put ${\mathcal O}_{E(H)}(-1) := {\tilde \pi}^*
(\wedge^t \tau)$, which is a line bundle 
on $E(H)$. 
Since $(\tau_{E(H)})^* \cong \wedge^{t-1}\tau_{E(H)} 
\otimes {\mathcal O}_{E(H)}(1)$, one can check that 
$$ \{\Sigma^{\alpha}
(\tau_{E(H)})^*\} = \{\Sigma^{\alpha}
(\tau_{E(H)})\}\otimes {\mathcal O}_{E(H)}(h-t) $$ 
where $\alpha$ runs through the Young diagrams 
with $r(\alpha) \leq t$ and $c(\alpha) \leq h-t$.  
The right hand side is a basis of $K(E(H))$.  
\vspace{0.12cm}

 {\bf Theorem (2.7)}. {\em The homomorphism 
 $$\Psi_0: K(T^*G^+) \to K(T^*G)$$ 
 is an isomorphism.}  
 \vspace{0.12cm}
 
 (2.7.1): 
 $E(H) \times_{\bar{\mathcal N}}E(H^+)$ is 
an integral scheme by [Ma 2, Corollary 3.15]. 
For example, when $t = 1$,  
 $E(H) \times_{\bar{\mathcal N}}E(H^+)$ coincides 
 with the blow-up of $E(H)$ ($E(H^*)$) along 
 the zero section.  
 When $t = 2$ and $h = 4$, 
$E(H) \times_{\bar{\mathcal N}}E(H^+)$ is a normal variety 
with only rational singularities (cf.(3.7)).
This implies that $E(H) \times_{\bar{\mathcal N}}E(H^+)$ 
is flat over the base space $\mathbf{C}^1$. 
We shall use this fact in the proof of (2.7).    
\vspace{0.12cm}
 
 (2.7.2).  Before proving (2.7) we review the 
 notion of a {\em restriction map with supports} 
 of Grothendieck groups (cf. [C-G, p. 246]). 
 Let $f: X_0 \to X$ be a closed immersion of 
 schemes with $X$ being a non-singular 
 quasi-projective variety.  
 For a closed subscheme $Z \subset X$, we 
 put $Z_0 : = Z \times_X X_0$. We shall  
 define the restriction map with the supports: 
 $$ f^*: K(Z) \to K(Z_0). $$ 
 Let $i: Z \to X$ be the inclusion map.      
 Let $F$ be a coherent sheaf on $Z$. 
 Since $i_*F$ has a finite locally free 
 resolution, $\mathrm{Tor}_n^{{\mathcal O}_X}
 ({\mathcal O}_{X_0}, i_*F)$ are zero except for a 
 finite number of $n$.  Moreover,  each 
 $\mathrm{Tor}_n^{{\mathcal O}_X}({\mathcal O}_{X_0}, i_*F)$ 
is an 
 ${\mathcal O}_{Z_0}$ module. Now we define 
 $f^*([F]) := \Sigma (-1)^n 
 [\mathrm{Tor}_n^{{\mathcal O}_X}
 ({\mathcal O}_{X_0}, i_*F)]$.   The restriction 
 map $f^*$ depends on the ambient spaces 
 $X$ and $X_0$.  
 \vspace{0.12cm}
 
 (2.7.3): {\em Proof of (2.7)}.  
 Consider the closed immersion 
 $$ f: T^*G \times T^*G^+ \to 
 E(H) \times_{{\mathbf C}^1}E(H^*). $$ 
 For $E(H) \times_{{\bar{\mathcal N}}}E(H^*) 
 \subset E(H) \times_{{\mathbf C}^1}E(H^*)$, 
 we define the restriction map  
 with supports:     
 $$ f^*: K(E(H)) \times_{{\bar{\mathcal N}}}E(H^*)) 
 \to K(T^*G \times_{\bar N}T^*G^+). $$ 
 On the other hand, one can define 
 the restriction map $i^*: K(E(H^+)) 
 \to K(T^*G^+)$ because $E(H^+)$ is 
 non-singular.  
 First let us check that the following diagram 
 commutes: 
 
 $$K(E(H^*)) \to K(T^*G^+)$$ 
 $$(\mu^+)^* \downarrow \hspace{0.6cm} 
 (\mu^+_0)^* \downarrow$$ 
 $$K(E(H) \times_{{\bar{\mathcal N}}}E(H^*)) 
 \to K(T^*G \times_{\bar N}T^*G^+). $$         
 
 It is sufficient to check the commutativity 
 for $[F] \in K(E(H^*))$ with $F$ a locally 
 free sheaf.  We have $(\mu^+_0)^* \circ 
 i^*([F]) = [(\mu^+_0)^*i^*F]$.  
 In turn, $f^* \circ (\mu^+)^*([F]) = 
 f^*([(\mu^+)^*F])$.  For simplicity, put  
 $X := E(H) \times_{{\mathbf C}^1}E(H^*)$,  
 $X_0 := T^*G \times T^*G^+$ 
 and $Z := E(H) \times_{{\bar{\mathcal N}}}E(H^*)$.  
 By definition,  $f^*([(\mu^+)^*F])
 = \Sigma (-1)^n 
 [\mathrm{Tor}_n^{{\mathcal O}_X}
 ({\mathcal O}_{X_0}, (\mu^+)^*F)]$. 
 Note that there is a surjective morphism 
 $X \to {\mathbf C}^1$ and its central fiber 
 is $X_0$. Since $X$ is flat over ${\mathbf C}^1$, 
 $\mathrm{Tor}_n^{{\mathcal O}_X}
 ({\mathcal O}_{X_0}, (\mu^+)^*F) =  
 \mathrm{Tor}_n^{{\mathcal O}_{{\mathbf C}^1}}
 (k(0), (\mu^+)^*F)$.  By (2.7.1), 
  $(\mu^+)^*F$ is a flat ${\mathcal O}_{{\mathbf C}^1}$ 
module.  
Therefore, $\mathrm{Tor}_n^{{\mathcal O}_{{\mathbf C}^1}}
 (k(0), (\mu^+)^*F) = 0$ when $n \ne 0$. 
 As a consequence,  we have 
 $f^*([(\mu^+)^*F])
 =  [((\mu^+)^*F)\vert_{X_0}]
 =  [(\mu^+_0)^*i^*F]$ 
 and the diagram commutes.     
 \vspace{0.12cm}
 
 Next, by [C-G, Proposition 5.3.15], the 
 following diagram commutes.  
 
 $$K(E(H) \times_{{\bar{\mathcal N}}}E(H^*)) 
 \to K(T^*G \times_{\bar N}T^*G^+) $$
 $$\mu_* \downarrow \hspace{0.6cm} 
 (\mu_0)_* \downarrow$$
 $$K(E(H)) \to K(T^*G).$$
 
 By the two commutative diagrams,  
 we see that homomorphisms $\Psi$ and 
 $\Psi_0$ are compatible with 
 the pull-backs:  
 
 $$K(E(H^*)) \to K(T^*G^+)$$ 
 $$\Psi \downarrow \hspace{0.6cm} 
 \Psi_0 \downarrow$$
 $$K(E(H)) \to K(T^*G).$$           
 
 By Corollary (2.2.3), the horizontal maps 
 are both isomorphisms. Since 
 $\Psi$ is an isomorphism by Theorem 
 (2.7),  $\Psi_0$ is also an isomorphism.  
 \vspace{0.12cm}

{\bf Example (2.8)}. Assume that $t = 1$ and 
$h = 3$. Then $G = {\mathbf P}^2$ 
and $G^+ = {\mathbf P}^2$. Note that $T^*G 
\times_{{\bar N}}T^*G^+$ is a normal crossing 
variety with 2 irreducible components. 
Let $X$ be the main component, namely, 
the irreducible component which dominates 
both $T^*G$ and $T^*G^+$. Let $p : 
X \to T^*G$ and $p^+: X \to T^*G^+$ 
be the projections. We define a homomorphism 
$$ {\Psi}'_0 : K(T^*G^+) \to K(T^*G) $$ 
as ${\Psi}'_0 := p_*\circ (p^+)^*$.     
In [Na], we have proved that, this ${\Psi}'_0$ 
is not an equivalence at the level of derived 
categories. Here we shall show that ${\Psi}'_0$ 
is not an isomorphism even at the level of K-theory.  

Put  ${\mathcal O}_{T^*G}(k) := \pi^*{\mathcal O}_G(k)$ 
and ${\mathcal O}_{T^*G^+}(k) := (\pi^+)^*{\mathcal O}_{G^+}(k)$.  
One can check that 
$$ {\Psi}'_0([{\mathcal O}_{T^*G^+}(-1)]) = 
[{\mathcal O}_{T^*G}(1)], $$ 
$$ {\Psi}'_0([{\mathcal O}_{T^*G^+}]) = 
[{\mathcal O}_{T^*G^+}], $$ 
and 
$$ {\Psi}'_0([{\mathcal O}_{T^*G^+}(1)]) = 
[{\mathcal O}_{T^*G}(-1)\otimes I].$$ 
Here $I$ is the ideal sheaf of the zero section of $\pi: 
T^*G \to G$. 
The Koszul resolution of $I$ yields the exact sequence 
$$ 0 \to \pi^*(\wedge^2\Theta_G \otimes 
{\mathcal O}_G(-1)) \to \pi^*(\Theta_G \otimes {\mathcal O}_G(-1)) 
\to I \otimes \pi^*{\mathcal O}_G(-1) \to 0. $$ 
By using the Euler exact sequence, we have 
$$ [\pi^*(\Theta_G \otimes {\mathcal O}_G(-1))] = 
3[{\mathcal O}_{T^*G}] - [{\mathcal O}_{T^*G}(-1)]$$ 
and 
$$ [\pi^*(\wedge^2\Theta_G \otimes 
{\mathcal O}_G(-1))] = 3[{\mathcal O}_{T^*G}(1)] 
- 3[{\mathcal O}_{T^*G}] 
+ [{\mathcal O}_{T^*G}(-1)]. $$
Therefore, we have 
$$ [I \otimes {\mathcal O}_{T^*G}(-1)] = 
-2[{\mathcal O}_{T^*G}(-1)] + 6[{\mathcal O}_{T^*G}] 
-3[{\mathcal O}_{T^*G}(1)]. $$  
The image of ${\Psi}'_0$ is the subgroup 
of $K(T^*G)$ generated by $[{\mathcal O}_{T^*G}(1)]$, 
$[{\mathcal O}_{T^*G}]$ and $2[{\mathcal O}_{T^*G}(-1)]$, 
which does not coincide with $K(T^*G)$.  
Hence, ${\Psi}'_0$ is {\em not} an isomorphism.  
\vspace{0.15cm}

\begin{center}
{\bf \S 3.} {\bf Stratified Mukai flop for G(2,4)}
\end{center}

Markman [Ma] has described a stratified Mukai flop 
as a sequence of blowing-ups and blowing-downs. 
Here we shall sketch this when $G = G(2, 4)$. 

(3.1). In the remainder of this section, we assume 
that $H$ is a 4-dimensional ${\mathbf C}$-vector space, 
$G = G(2, H)$ and $G^+ = G(2, H^*)$. For $p \in G$, 
let $\mathrm{Hom}^i(q(p), \tau(p)) := \{\phi \in 
\mathrm{Hom}(q(p), \tau(p)); \mathrm{rank}\phi 
\leq i \}$ and put 
${\underline{\mathrm{Hom}}}^i(q, \tau) := 
\cup_{p \in G}\mathrm{Hom}^i(q(p), \tau(p))$. 
By (1.3) we have a sequence of subvarieties of 
$E(H)$: 
$$ {\underline{\mathrm{Hom}}}^0(q, \tau) \subset 
{\underline{\mathrm{Hom}}}^1(q, \tau) \subset T^*G \subset 
E(H).$$   
For short, we write ${\mathcal X}$ for $E(H)$, 
$X$ for $T^*G$, $Z$ for ${\underline{\mathrm{Hom}}}^1(q, \tau)$ 
and $M$ for ${\underline{\mathrm{Hom}}}^0(q, \tau)$. 
As in (2.2), ${\tilde \pi}: {\mathcal X} \to G$ and 
$\pi: X \to G$ are projection maps. Note that 
$M$ is the zero section of 
${\tilde \pi}: {\mathcal X} \to G$.  
We have $\dim {\mathcal X} = 9$, $\dim X = 8$, 
$\dim Z = 7$ and $\dim M = 4$.   
Let ${\tilde s}: {\mathcal X} \to {\bar{\mathcal N}}^2(H)$ 
be the extended Steinberg resolution (cf. (1.5)). 
Let $\Sigma$ be the singular locus of 
${\bar{\mathcal N}}^2(H)$. Every element of 
$\Sigma$ is conjugate to  
a matrix 

$$ A = \left(
\begin{array}{cc}
\alpha I & *\\
0 & 0  
\end{array}\right)$$ 
whose rank $\leq 1$ (cf. (1.4)). 
This implies that $\alpha = 0$; hence  
$A^2 = 0$. Thus, we have    
$\Sigma \cong {\bar N}^1(H)$. Note 
that $\dim \Sigma = 6$.
The exceptional locus ${\mathrm 
Exc}({\tilde s})$ coincides with 
$Z$. There is a fibration 
${\tilde s}\vert_Z: 
Z \to \Sigma$. For $p \in 
\Sigma\setminus\{0\}$, ${\tilde s}^{-1}(p) 
\cong {\mathbf P}^1$. 
Clearly, ${\tilde s}^{-1}(0) = M$.         
\vspace{0.12cm}

(3.2).   
Let $\nu_1: {\mathcal X}_1 \to {\mathcal X}$ 
be the blowing up of ${\mathcal X}$ along 
$M$. 

(3.2.1). Let $\mathrm{Exc}(\nu_1)$ be 
the exceptional locus of $\nu_1$. Then 
$\mathrm{Exc}(\nu_1) \cong 
{\mathbf P}(E(H))$, which 
is a ${\mathbf P}^4$ bundle over $M$. 

(3.2.2). Let ${\tilde Z}$ be the proper 
transform of $Z$ by $\nu_1$. One can check 
that $Z$ has 3-dimensional ordinary double 
points along $M$. Hence $(\nu_1)\vert_{\tilde Z}: 
{\tilde Z} \to Z$ is a resolution of $Z$. 
The exceptional locus 
$\mathrm{Exc}((\nu_1)\vert_{\tilde Z})$ is 
a ${\mathbf P}^1 \times {\mathbf P}^1$ bundle 
over $M$.    

(3.2.3).  
${\tilde Z} \cap \mathrm{Exc}(\nu_1)$ is 
described as 
follows. Let $p \in M (\cong G)$. If we choose a 
suitable basis of $H$, then 
$E(H)(p)$ consists of the matrices of the 
following form  

\[\left(
\begin{array}{cccc}
\alpha & 0 & x & y \\
0 & \alpha & z & w \\ 
0 & 0 & 0 & 0 \\ 
0 & 0 & 0 & 0   
\end{array}\right)\] 

We can regard $(\alpha:x:y:z:w)$ 
as homogeneous coordinates of 
the projective space 
${\mathbf P}(E(H)(p))$. 
Then 
$${\tilde Z} \cap {\mathbf P}(E(H)(p))  
= \{\alpha = xw - yz = 0 \},$$ which 
is isomorphic to ${\mathbf P}^1 
\times {\mathbf P}^1$. In this way, 
${\tilde Z} \cap {\mathbf P}(E(H))$ 
becomes a 
${\mathbf P}^1 \times {\mathbf P}^1$ 
bundle over $M$. 

(3.2.4). Let $T^*G(1,H) \to 
{\bar N}^1(H)$ be the Springer 
resolution (cf. (1.2)). Let ${\tilde N}$ 
be the blowing up of $T^*G(1,H)$ along 
the zero section. Identify ${\bar N}^1(H)$ 
with $\Sigma$ (cf.(3.1)). Note that 
${\tilde Z} \to Z$ is the blowing-up 
along $M$, and ${\tilde N} \to \Sigma$ 
is the blowing-up at $0 \in \Sigma$. 
Then the composite ${\tilde Z} \to 
Z \to \Sigma$ is factorized as
 
$$ {\tilde Z} \to Z $$ 
$$ \downarrow \hspace{0.7cm} \downarrow $$ 
$$ {\tilde N} \to \Sigma. $$ 

Moreover, ${\tilde Z} \to {\tilde N}$ 
is a ${\mathbf P}^1$ bundle. The proof 
of these facts are omitted.  
\vspace{0.12cm}

(3.3). Let $\nu_2: {\tilde{\mathcal X}} \to 
{\mathcal X}_1$ be the blowing up of 
${\mathcal X}_1$ along ${\tilde Z}$. 
We put $E := \mathrm{Exc}(\nu_2)$. 
$E$ is a ${\mathbf P}^1$ bundle over 
${\tilde Z}$. 
Moreover, let $F$ be the proper 
transform (=total transform) of 
$\mathrm{Exc}(\nu_1)$. 
$(\nu_1\circ\nu_2)\vert_F: F \to 
M$ is a smooth morphism, whose 
fibers are blown-up ${\mathbf P}^4$ 
along ${\mathbf P}^1 \times {\mathbf P}^1$ 
described in (3.2.3).  
\vspace{0.12cm}

(3.4). For the dual $G^+ = G(2,H^*)$, 
we also have the varieties 
$$ M^+ \subset Z^+ \subset X^+ \subset {\mathcal X}^+.$$ 
By the same way as (3.1),(3.2) and (3.3), we 
have a sequence of blowing-ups: 
$$ {\tilde{\mathcal X}}^+ \stackrel{\nu^+_2}\to 
{\mathcal X}^+_1 \stackrel{\nu^+_1}\to {\mathcal X}^+.$$ 
The birational map ${\mathcal X} - - \to {\mathcal X}^+$ 
in (1.6) 
induces a birational map 
${\tilde{\mathcal X}} - - \to {\tilde{\mathcal X}}^+$.  
This birational map is actually an isomorphism  
by [Ma]. \vspace{0.12cm}       

(3.5). We put ${\bar{\mathcal N}} := 
{\bar{\mathcal N}}^2(H) \cong 
{\bar{\mathcal N}}^2(H^*)$. 
The fiber product 
${\mathcal X}\times_{{\bar{\mathcal N}}}
{\mathcal X}^+$ birationally 
dominates ${\mathcal X}$.  
Now we have a natural birational 
map over ${\mathcal X}$: 
$${\mathcal X}_1 - - \to 
{\mathcal X}\times_{{\bar{\mathcal N}}}
{\mathcal X}^+.$$ 
Composing this with $\mu_2: 
{\tilde{\mathcal X}} \to {\mathcal X}_1$, 
we have a birational map 
$${\tilde{\mathcal X}} - - \to 
{\mathcal X}\times_{{\bar{\mathcal N}}}
{\mathcal X}^+,$$ 
which is, in fact, a morphism 
because ${\tilde{\mathcal X}} - - \to 
{\mathcal X}^+$ is a morphism.
Note that 
these birational maps are 
defined over ${\bar{\mathcal N}}$. 
Over $0 \in {\bar{\mathcal N}}$, 
the first birational map 
induces a rational map 
$$ \gamma : {\mathbf P}(E(H)) 
- - \to M \times M^+, $$ 
and the second birational morphism  
induces a morphism 
$$ F \rightarrow M \times M^+.$$ 
We shall describe this last map 
as the blowing-up 
of $M \times M^+$ with a suitable 
center.  
\vspace{0.12cm}
   
(3.5.1). Let  ${\gamma}(p): 
{\mathbf P}(E(H)(p)) - - \to M^+$ be 
the restriction of $\gamma$ 
over $p \in M$. 
We shall give an explicit description of 
 $\gamma(p)$.  By a suitable choice of 
the basis of $H$,  $E(H)(p)$ is 
identified with the set of matrices in (3.2.3). 
Here we introduce the dual basis in 
$H^*$, which will be used later.  
Fix a matrix 

$$\left(
\begin{array}{cccc}
\alpha & 0 & x & y \\ 
0 & \alpha & z & w \\ 
0 & 0 & 0 & 0 \\ 
0 & 0 & 0 & 0 
\end{array}\right)$$ 
in $E(H)(p)$ in such a way 
that  $\alpha \ne 0$. 
Let $l := \{q(t)\}_{t 
\in {\mathbf C}}$ be the line 
of $E(H)(p)$ 
passing through this matrix and 
$0$. Here    

$$q(t) :=  
\left(
\begin{array}{cccc}
\alpha t & 0 & xt & yt \\ 
0 & \alpha t & zt & wt \\ 
0 & 0 & 0 & 0 \\ 
0 & 0 & 0 & 0 
\end{array}\right)$$ 
is a 1-parameter family of 
matrices.     
Let $[l] \in {\mathbf P}(E(H)(p))$ 
be the point represented by 
$l$. The image $\gamma(p)([l]) \in 
M^+$ of $[l]$ is described as follows. 
We have the following commutative 
diagram: 

$$ E(H)  - - \to E(H^*) $$ 
$$ {\tilde s}\downarrow \hspace{0.5cm} 
{\tilde s}^+ \downarrow $$ 
$$ {\bar{\mathcal N}}^t(H) \cong 
{\bar{\mathcal N}}^t(H^*).$$
 
For $(p, q(t)) \in E(H)$, ${\tilde s}((p, q(t)) 
= q(t)$.  The isomorphism  
${\bar{\mathcal N}}^2(H) \cong 
{\bar{\mathcal N}}^2(H^*)$ is given 
by the transposition. Hence this 
isomorphism sends $q(t)$ to 

$${ }^tq(t) := 
\left(
\begin{array}{cccc}
\alpha t & 0 & 0 & 0 \\ 
0 & \alpha t & 0 & 0 \\ 
xt & zt & 0 & 0 \\
yt & wt & 0 & 0 
\end{array}\right).$$ 
Assume that $t \ne 0$. 
Then, since $\alpha \ne 0$, 
the inverse image $({\tilde s}^+)^{-1}({ }^tq(t))$ 
is uniquely determined and is given by 
$(?(t), { }^tq(t)) \in E(H^*)$. Here 
$?(t) \in G^+$ is the 2-dimensional 
subspace of $H^*$ generated by 
two vectors 
 
$$\left(
\begin{array}{c}
  \alpha t \\ 
  0  \\ 
  xt \\ 
  yt  
\end{array}\right), 
\left(
\begin{array}{c} 
  0 \\ 
  \alpha t \\ 
  zt  \\ 
  wt 
\end{array}\right).$$

Now ${\mathrm lim}_{t \to 0} 
(?(t), { }^tq(t)) = (\gamma(p)([l]), 0)$. 
Therefore, $\gamma(p)([l]) \in M^+(\cong G^+)$ 
is the 2-dimensional subspace of 
$H^*$ generated by two vectors 

$$\left(
\begin{array}{c}
  \alpha \\ 
  0  \\ 
  x \\ 
  y
\end{array}\right),    
\left(
\begin{array}{c} 
  0 \\ 
  \alpha \\ 
  z  \\ 
  w 
\end{array}\right).$$

Let $M^+  \to {\mathbf P}^5$ 
be the Pl{\"u}cker embedding with 
respect to the basis of $H^*$ given 
above. Let $(p_{12}: p_{13}: p_{14}: 
p_{23}: p_{24}: p_{34})$ be the Pl{\"u}cker 
coordinates. 
Note that $M^+$ is a quadratic hypersurface 
of ${\mathbf P}^5$ defined by 
$p_{12}p_{34} - p_{13}p_{24} + p_{14}p_{23} = 0$. 
Then $\gamma(p): 
{\mathbf P}(E(H)(p))  - - \to M^+ \subset 
{\mathbf P}^5$ is given by 
$$ \gamma(p)((\alpha : x : y : z : w)) = 
({\alpha}^2: \alpha z: \alpha w : -\alpha x : 
-\alpha y : xw - yz). $$    
The indeterminancy of $\gamma(p)$ is 
the subvariety 
$\{{\alpha} = x w - yz = 0\}$. This 
coincides with  ${\tilde Z} \cap 
{\mathbf P}(E(H)(p))$ (cf.  (3.2.3)). 
Let $F(p)$ be the fiber of the 
morphism $F \to M$ over $p \in M$ 
(cf. (3.3)). Then,  
$F(p)$ is the blow-up of ${\mathbf P}
(E(H)(p))$ along ${\tilde Z} \cap 
{\mathbf P}(E(H)(p))$.  It is immediately 
checked that the composite 
$F(p) \to {\mathbf P}(E(H)(p)) 
- - \to M^+$ is a birational morphism.    
The following are also checked.
\vspace{0.12cm}

(3.5.1-a).  Let $R(p)$ be the proper transform of $\{\alpha = 0\} 
\subset {\mathbf P}(E(H)(p))$ by the blowing-up 
$F(p) \to {\mathbf P}(E(H)(p))$.  Then $R(p)$ is isomorphic to 
$\{\alpha = 0\}(\cong {\mathbf P}^3)$. Moreover, 
$R(p)$ is contracted to the point $(0:0:0:0:0:1) \in 
{\mathbf P}^5$ by the birational morphism 
$F(p) \to M^+$.  Actually, $F(p)$ is the blowing-up 
of $M^+$ at $(0:0:0:0:0:1)$.   \vspace{0.12cm}

(3.5.1-b).  Let $S(p)$ be the exceptional divisor of the 
blowing-up $F(p) \to {\mathbf P}(E(H)(p))$.  Then 
$S(p)$ is mapped onto the divisor $\{p_{12} = 0\} 
\cap M^+$ of $M^+$ by the birational map 
$F(p) \to M^+$. This divisor has an ordinary 
double threefold singularity at $(0:0:0:0:0:1)$.  
Actually, $S(p)$ is the blowing-up 
of $\{p_{12} = 0\} \cap M^+$ at $(0:0:0:0:0:1)$. 
The exceptional locus of this blowing-up is 
$S(p) \cap R(p)$, which is isomorphic to 
${\mathbf P}^1 \times {\mathbf P}^1$.     
\vspace{0.12cm}

(3.5.2). In general, we have a 
natural isomorphism $G(t, H) 
\cong G(h-t, H^*)$, where $\dim 
H = h$. Now, since $t = 2$ and $h = 4$, 
there is an isomorphism 
$\iota: G(2, H) \cong G(2, H^*)$.  
The center $p^+ \in M^+$ of the blowing-up  
$F(p) \to M^+$ depends on $p \in M$. 
By (3.5.1-a), we see that $p^+ = 
\iota(p)$. Let $\Gamma \subset 
M \times M^+$ be the graph of 
$\iota$. By (3.5.1) the birational map 
$F \to M \times M^+$ is the blowing-up 
of $M \times M^+$ along $\Gamma$: 
$$ F \cong 
{\mathrm Bl}_{\Gamma}(M \times M^+). $$ 

(3.6).  As in (3.4.1), for ${\tilde{\mathcal X}}^+$, 
there is a birational morphism 
$${\tilde{\mathcal X}}^+  \rightarrow
{\mathcal X} \times_{{\bar{\mathcal N}}}{\mathcal X}^+ $$ 
as ${\bar{\mathcal N}}$ schemes. Over 
$0 \in {\bar{\mathcal N}}$, this birational 
morphism induces a morphism 
$$ \gamma^+ : F^+  \rightarrow 
M \times M^+. $$  
By the same argument as (3.5.2) we see that 
$$ F^+ \cong  
{\mathrm Bl}_{\Gamma}(M \times M^+). $$ 

(3.7). One can 
check that 
${\mathcal X} \times_{{\bar{\mathcal N}}}{\mathcal X}^+$ 
is a normal variety with rational singularities. 
In fact, set-theoretically, 
${\mathcal X} \times_{{\bar{\mathcal N}}}{\mathcal X}^+$ 
is obtained from ${\tilde{\mathcal X}}$ by 
contracting $\{R(p)\}_{p \in M}$ to $\{p \times 
p^+ \in M \times M^+\}_{p \in M}$. Therefore, 
${\mathcal X} \times_{{\bar{\mathcal N}}}{\mathcal X}^+$ 
has singularities along $\Gamma$ in (3.5.2). 
Moreover, by a direct calculation, we see that 
these singularities are locally of the following 
form:  
$(V, 0) \times ({\mathbf C}^4, 0),$ 
where 
$$V := \{(x,y,z,w,s,t,u,v) \in {\mathbf C}^8; 
\mathrm{rank}
\left(
\begin{array}{cccc}
x & y & z & w \\ 
-v & t & u & -s 
\end{array}\right) \leq 1\}.$$  
\vspace{0.12cm}

(3.8).  Let $\Sigma$ be the singular locus of 
${\bar{\mathcal N}}$ and put $\Sigma^* := 
\Sigma \setminus \{0\}$. $Z$ (resp. $Z^+$) 
has a fibration over $\Sigma$ (cf. (3.1)). Let 
$Z^*$ (resp. $(Z^+)^*$) be the 
inverse image of $\Sigma^*$ 
by this map. By (3.1), $Z^* \to 
\Sigma^*$ and $(Z^+)^* \to 
\Sigma^*$ are both ${\mathbf P}^1$ 
bundles. Outside $M$ and $M^+$, the 
birational map ${\mathcal X} - - \to {\mathcal X}^+$ is  
 a family of Atiyah flops along $Z^* \to 
\Sigma^*$ and $(Z^+)^* \to \Sigma^*$. 
\vspace{0.15cm} 

\begin{center}
{\bf \S 4. Derived Categories}
\end{center}

Let 
$$ \Psi : D(E(H^*)) \to D(E(H)) $$ 
be the functor defined by $\Psi(\bullet) 
:= {\mathbf R}\mu_*\circ 
{\mathbf L}(\mu^+)^*(\bullet)$(cf. (2.4)), 
where 
$D(E(H))$ (resp. $D(E(H^*))$) is 
the bounded derived category of 
coherent sheaves on $E(H)$
(resp. $E(H^*)$). 
Let 
$$ \Psi_0 : D(T^*G^+) 
\to D(T^*G) $$ 
be the functor defined by $\Psi_0(\bullet) 
:= {\mathbf R}(\mu_0)_*\circ 
{\mathbf L}(\mu_0^+)^*(\bullet)$(cf. (2.4)).  
In this section, we show that these 
functors are not equivalences when 
$G = G(2,4)$.  
In the remainder, we use the same 
notation as Section 3.  
\vspace{0.12cm}

{\bf Lemma (4.1)}. {\em Let 
$pr: F \to M$ and $pr^+: F \to M^+$ 
be the projections in (3.5). 
Then 
$$ (pr^+)^*{\mathcal O}_{M^+}(1) 
\otimes {\mathcal O}_F(E + 2F) 
\in pr^*\mathrm{Pic}(M), $$ 
where ${\mathcal O}_{M^+}(1)$ 
is the tautological line bundle 
of the Grassmannian $M^+$.}  
\vspace{0.12cm}

{\em Proof}. It is sufficient to prove 
that, for $p \in M$, 
$$(pr^+)^*{\mathcal O}_{M^+}(1) 
\otimes {\mathcal O}_{F(p)}(E + 2F) 
\cong {\mathcal O}_{F(p)}. $$ 
Here $F(p)$ is the 
blowing-up of ${\mathbf P}^4(=  
{\mathbf P}(E(H)(p)))$ along 
${\tilde Z} \cap 
{\mathbf P}(E(H)(p))$. We call 
this blowing-up $\nu_2 (p)$ 
(3.5.1).  Note that ${\mathrm 
Exc}(\nu_2 (p)) = E \cap F(p)$. 
Then, by (3.5.1-a) and (3.5.1-b), 
$$(pr^+\vert_{F(p)})^*{\mathcal O}_{M^+}(1) = 
(\nu_2)^*{\mathcal O}_{{\mathbf P}^4}(2) \otimes 
{\mathcal O}_{F(p)}(-E) = {\mathcal O}_{F(p)}(-2F-E).$$  
\vspace{0.12cm}

(4.2) Recall that $X := T^*G$ and 
$\pi: X \to G$ is the projection. 
Since $X \cong 
{\underline{\mathrm{Hom}}}(q, \tau)$, 
there is a universal homomorphism 
$$f_{univ}: \pi^*q \to \pi^*\tau.$$ 
Then $Z$ is the divisor of $X$ defined 
by $\wedge^2 f_{univ} = 0$, where  
$\wedge^2f_{univ}$ is an element of 
$\Gamma (X, {\underline{\mathrm{Hom}}}
(\wedge^2\pi^*q, \wedge^2\pi^*\tau)) = 
\Gamma (X, \pi^*{\mathcal O}_G(-2)).$ 
We have an exact sequence: 
$$ 0 \to \pi^*{\mathcal O}_G(2) \to 
{\mathcal O}_X \to {\mathcal O}_Z \to 0.$$ 
By (1.5), $X$ is the central fiber of 
the morphism 
${\mathcal X}:= E(H) \to 
{\mathbf C}^1$. 
Hence, there is an exact 
sequence: 
$$ 0 \to {\mathcal O}_{{\mathcal X}} \to 
{\mathcal O}_{{\mathcal X}} \to {\mathcal O}_X \to 0.$$  @@@ 
   
(4.3). Put $\nu := \nu_1\circ\nu_2$ 
and $\nu^+ := (\nu^+_1)\circ(\nu^+_2)$. 
Then we have a diagram 
$$ {\mathcal X} \stackrel{\nu}\leftarrow 
{\tilde{\mathcal X}} \cong 
{\tilde{\mathcal X}^+} \stackrel{\nu^+}
\rightarrow {\mathcal X}^+.$$ 
This diagram is defined over 
the parameter space ${\mathbf C}^1$. 
The restriction of this to the 
central fibers becomes 
$$ X \stackrel{\nu_0}\leftarrow 
{\tilde X} \cong 
{\tilde X}^+ \stackrel{\nu_0^+}
\rightarrow X^+.$$
By the isomorphism 
${\tilde{\mathcal X}} \cong 
{\tilde{\mathcal X}^+}$, we regard 
$\nu^+$ as a morphism from 
${\tilde{\mathcal X}}$ to 
${\mathcal X}^+$. We simply write 
${\mathcal O}(k)$ for ${\tilde \pi}^*
{\mathcal O}_G(k)$, and ${\mathcal O}^+(k)$ 
for $({\tilde \pi}^+)^*{\mathcal O}_{G^+}(k)$. 
Define a functor
$\Phi$ as ${\mathbf R}{\nu}_*\circ 
{\mathbf L}{\nu^+}^*$ and 
a functor $\Phi_0$ as 
${\mathbf R}{\nu_0}_*\circ 
{\mathbf L}{\nu_0^+}^*$:  
$$ \Phi: D({\mathcal X}^+) \to 
D({\mathcal X}), $$ 
$$ \Phi_0: D(X^+) \to 
D(X).$$   

(4.4). We first show that 
$\Phi$ is not an equivalence. 
Since ${\bar{\mathcal N}}^2(H)$ is a 
Stein space with rational singularities 
and ${\mathcal X}$ is a resolution of 
${\bar{\mathcal N}}^2(H)$(cf.(1.5)), we have  
$\mathrm{Ext}^i({\mathcal O}^+(1), 
{\mathcal O}^+(1)) = 0$ for $i > 0$.  
Let us compute 
$\mathrm{Ext}^i(\Phi({\mathcal O}^+(1)), 
\Phi({\mathcal O}^+(1))).$  
\vspace{0.12cm}

(4.5). By (4.1) and (3.8), we see 
that the restriction of 
$(\nu^+)^*{\mathcal O}^+(1)(E+2F)$ to 
each fiber of $\nu$ is a trivial 
line bundle. This implies that 
$(\nu^+)^*{\mathcal O}^+(1)(E+2F)$ is 
a trivial line bundle around each 
fiber of $\nu$. Hence, $\nu_*
(\nu^+)^*{\mathcal O}^+(1)(E+2F)$ is 
a line bundle on ${\mathcal X}$ and 
its pull-back by $\nu$ coincides 
with $(\nu^+)^*{\mathcal O}^+(1)(E+2F)$. 
By (2.3.2), $\nu_*
(\nu^+)^*{\mathcal O}^+(1)(E+2F) \cong 
{\mathcal O}(-1)$. Hence we have 
$(\nu^+)^*{\mathcal O}^+(1)(E+2F) \cong 
\nu^*{\mathcal O}(-1).$
Now the exact sequence 
$$ 0 \to (\nu^+)^*{\mathcal O}^+(1) 
\to (\nu^+)^*{\mathcal O}^+(1)(E+2F) 
\to (\nu^+)^*{\mathcal O}^+(1)\vert_{E+2F} 
\to 0 $$ 
is identified with the exact sequence 
$$ 0 \to \nu^*{\mathcal O}(-1)(-E-2F) 
\to \nu^*{\mathcal O}(-1) \to 
\nu^*{\mathcal O}(-1)\vert_{E+2F} \to 0.$$ 
Apply $\nu_*$ to this last sequence. 
Since $R^1\nu_*(\nu^+)^*{\mathcal O}^+(1) 
= 0$, we have the exact sequence 
$$ 0 \to {\mathcal O}(-1)\otimes 
\nu_*{\mathcal O}(-E-2F) \to 
{\mathcal O}(-1) \to {\mathcal O}(-1)
\otimes \nu_*{\mathcal O}_{E+2F} \to 
0.$$ Let $Z'$ be the scheme theoretic 
image of $E+2F$ by $\nu$. By definition, 
the ideal sheaf $I_{Z'}$ of $Z'$ is 
$\nu_*{\mathcal O}(-E-2F)$. 
Then, this sequence is obtained 
from the exact sequence 
$$ 0 \to I_{Z'} \to {\mathcal O}_{\mathcal X} \to 
{\mathcal O}_{Z'} \to 0 $$ 
by taking the tensor product with 
${\mathcal O}(-1)$. 
Now 
$Z$ coincides with the scheme theoretic image of 
$E+F$. By the next lemma, 
we have an exact sequence 

$(4.5.1): \hspace{1.5cm} 0 \to {\mathcal O}_M  
\to {\mathcal O}_{Z'} \to 
{\mathcal O}_Z \to 0.$ 
\vspace{0.12cm}

Since $R^i\nu_*(\nu^+)^*{\mathcal O}^+(1) 
= 0$ for $i > 0$, we have 
$$ \Phi({\mathcal O}^+(1)) = 
{\mathcal O}(-1) \otimes I_{Z'}.$$ 
\vspace{0.12cm}

{\bf Lemma (4.5.2)}. {\em There  
is an exact sequence 
$$ 0 \to \nu_*{\mathcal O}_F(-E-F) 
\to {\mathcal O}_{Z'} \to {\mathcal O}_Z 
\to 0,$$ and $\nu_*{\mathcal O}_F(-E-F) 
\cong {\mathcal O}_M.$} 
\vspace{0.12cm}

{\em Proof}. The first claim easily follows 
from the definitions of $Z$ and $Z'$.  
For the second claim, we first show that 
$\nu_*{\mathcal O}_F(-E-F)$ is a line bundle 
on $M$. It is enough to prove that, for 
$p \in M$, $h^0(F(p), {\mathcal O}_{F(p)}(-E-F)) 
= 1.$ By (3.5.1), $F(p)$ is the blowing-up 
of ${\mathbf P}(E(H)(p))$ along $e(p):= 
{\tilde Z} \cap {\mathbf P}(E(H)(p))$. We 
call this blowing-up $\nu_2(p)$ as in the 
proof of (4.1).  
We shall use the homogenous coordinates 
$(\alpha:x:y:z:w)$ of ${\mathbf P}^4 = 
{\mathbf P}(E(H)(p))$ 
in (3.5.1). Then $e(p)= \{\alpha = xw-yz = 0\}.$ 
Let $I_{e(p)}$ be the ideal sheaf of 
$e(p)$ in ${\mathbf P}^4$. 
Since ${\mathcal O}_{F(p)}(F) = (\nu_2(p))^*
{\mathcal O}_{{\mathbf P}^4}(-1)$, 
we only have to prove that $h^0({\mathbf P}^4, 
{\mathcal O}_{{\mathbf P}^4}(1)\otimes I_{e(p)}) 
= 1$. But this is checked directly.  
We next show that 
$\nu_*{\mathcal O}_F(-E-F)$ has a nowhere-vanishing 
section. Let $R(p)$ be the same as (3.5.1-a). 
Then $R(p)$ is a non-zero section of 
$H^0((F(p), {\mathcal O}_{F(p)}(-E-F))$. Now, 
$R := \{R(p)\}_{p \in M}$ gives a 
nowhere-vanishing section of  
$\nu_*{\mathcal O}_F(-E-F)$.   
\vspace{0.12cm}

(4.6). By (4.5) we have 
$$\mathrm{Ext}^i(\Phi({\mathcal O}^+(1)), 
\Phi({\mathcal O}^+(1))) \cong 
 \mathrm{Ext}^i(I_{Z'}, I_{Z'}).$$ 
We shall prove that 
$\mathrm{Ext}^5(I_{Z'}, I_{Z'}) \ne 0.$ 
\vspace{0.12cm}

{\bf Lemma (4.6.1)}. {\em $H^i({\mathcal X}, I_{Z'}) 
 = 0$ for $i > 1$.}  
\vspace{0.12cm}

{\em Proof}. Let $I_Z$ be the ideal sheaf 
of $Z$. Then, by (4.5.1), there is an exact 
sequence 
$$ 0 \to I_{Z'} \to I_Z \to {\mathcal O}_M \to 
0.$$ 
Since $H^i({\mathcal O}_M) =  
 0$ for $i > 0$, it is enough to prove that 
$H^i(I_Z) = 0$ for $i > 0$. We use the diagram 
in (3.2.4):   
  
$$ {\tilde Z} \to Z $$ 
$$ \downarrow \hspace{0.4cm} \downarrow $$ 
$$ {\tilde N} \to \Sigma. $$ 

Since $Z$ has rational singularities 
(cf.(3.2.2)) and $\Sigma \cong {\bar N}^1(H)$ also 
has a rational singularity, we have 
$H^i({\mathcal O}_Z) = 0$ for $i > 0$. 
Then the results follow from the 
exact sequence 
$$ 0 \to I_Z \to {\mathcal O}_{\mathcal X} 
\to {\mathcal O}_Z \to 0 $$ 
because $H^i({\mathcal O}_{\mathcal X}) = 0$ 
for $i > 0$.  
\vspace{0.12cm}

{\bf Lemma (4.6.2)}. 
$$\mathrm{Ext}^i({\mathcal O}_X, 
{\mathcal O}_{\mathcal X}) = \left\{ 
\begin{array}{rl}    
0 & \quad (i \ne 1)\\ 
H^0(X, {\mathcal O}_X) & \quad 
(i = 1)
\end{array}\right.$$ 

$$\mathrm{Ext}^i({\mathcal O}_X, 
{\mathcal O}(-2)) = \left\{ 
\begin{array}{rl}    
0 & \quad (i \ne 1,2)\\ 
H^0(X, {\mathcal O}(-2)\vert_X) & \quad 
(i = 1) \\
H^1(X, {\mathcal O}(-2)\vert_X) & \quad 
(i = 2)
\end{array}\right.$$ 
   
{\em Proof}. By the exact sequence in (4.2): 
$$ 0 \to {\mathcal O}_{\mathcal X} \to 
{\mathcal O}_{\mathcal X} \to {\mathcal O}_X \to 0,$$ 
we have 

$${\underline{\mathrm{Ext}}}^i({\mathcal O}_X, 
{\mathcal O}_{\mathcal X}) = \left\{ 
\begin{array}{rl}    
0 & \quad (i \ne 0)\\ 
{\mathcal O}_X & \quad 
(i = 1)
\end{array}\right.$$ 

$${\underline{\mathrm{Ext}}}^i({\mathcal O}_X, 
{\mathcal O}(-2)) = \left\{ 
\begin{array}{rl}    
0 & \quad (i \ne 1)\\ 
{\mathcal O}(-2)\vert_X & \quad 
(i = 1)
\end{array}\right.$$ 

Note that $H^j({\mathcal O}(-2)\vert_X) 
= 0$ for $j \geq 2)$, and 
$H^j({\mathcal O}_X) = 0$ for 
$j \geq 1$. The results follow from 
the local to global spectral sequence 
of Ext. \vspace{0.12cm}

{\bf Lemma (4.6.3)}. 
$$\mathrm{Ext}^i({\mathcal O}_Z, 
{\mathcal O}_{\mathcal X}) = \left\{ 
\begin{array}{rl}    
0 & \quad (i \ne 3)\\ 
H^1(X, {\mathcal O}(-2)\vert_X) & \quad 
(i = 3)
\end{array}\right.$$ 
         
{\em Proof}. Apply $\mathrm{Ext}(\bullet, 
{\mathcal O}_{\mathcal X})$ to the exact sequence 
$$ 0 \to {\mathcal O}(2)\vert_X \to {\mathcal O}_X 
\to {\mathcal O}_Z \to 0.$$ 
Since $\mathrm{Ext}^2({\mathcal O}_X, 
{\mathcal O}_{\mathcal X}) = 0$ by (4.6.2), we 
have an exact sequence 
$$ \mathrm{Ext}^1({\mathcal O}_X, 
{\mathcal O}_{\mathcal X}) \to 
\mathrm{Ext}^1({\mathcal O}(2)\vert_X, 
{\mathcal O}_{\mathcal X}) \to 
\mathrm{Ext}^2({\mathcal O}_Z, 
{\mathcal O}_{\mathcal X}) \to 0. $$ 
By (4.6.2) and (4.2), this sequence 
is identified with the exact sequence 
$$ H^0(X, {\mathcal O}_X) \to 
H^0(X, {\mathcal O}(-2)\vert_X) 
\to H^0(Z, {\mathcal O}(-2)\vert_Z) 
\to 0.$$ 

It is easily checked that 
$H^0(Z,{\mathcal O}(-2)\vert_Z) = 0$. 
This implies that  
$\mathrm{Ext}^2({\mathcal O}_Z,
{\mathcal O}_{{\mathcal X}}) = 0.$  
Since $\mathrm{Ext}^2({\mathcal O}_X, 
{\mathcal O}_{{\mathcal X}}) 
= \mathrm{Ext}^3({\mathcal O}_X,
{\mathcal O}_{{\mathcal X}}) = 0$, we see 
that 
$\mathrm{Ext}^2({\mathcal O}(2)\vert_X, 
{\mathcal O}_{{\mathcal X}}) \cong  
\mathrm{Ext}^3({\mathcal O}_Z, 
{\mathcal O}_{{\mathcal X}}).$  
By (4.6.2), 
$\mathrm{Ext}^2({\mathcal O}(2)\vert_X, 
{\mathcal O}_{{\mathcal X}}) = 
H^1(X, {\mathcal O}(-2)\vert_X).$
\vspace{0.12cm}

{\bf Lemma (4.6.4)}. 
$$\mathrm{Ext}^i({\mathcal O}_Z, 
{\mathcal O}_Z) = \left\{ 
\begin{array}{rl}    
0 & \quad (i \ne 0,3)\\ 
H^0(Z, {\mathcal O}_Z) & \quad 
(i = 0) \\
H^1(Z, {\mathcal O}(-2)\vert_Z) & \quad 
(i = 3)
\end{array}\right.$$ 

{\em Proof}. Since $Z \subset 
{\mathcal X}$ is locally of complete 
intersection, ${\underline{\mathrm 
Ext}}^i({\mathcal O}_Z, {\mathcal O}_Z)  
\cong \wedge^i N_{Z/{\mathcal X}}$. 
Since $N_{Z/X} \cong 
{\mathcal O}(-2)\vert_Z$ and 
$N_{X/{\mathcal X}}\vert_Z \cong 
{\mathcal O}_Z$, there is an exact 
sequence 
$$ 0 \to {\mathcal O}(-2)\vert_Z 
\to N_{Z/{\mathcal X}} \to {\mathcal O}_Z 
\to 0.$$ 
From this sequence we know that 
$${\underline{\mathrm{Ext}}}^i({\mathcal O}_Z, 
{\mathcal O}_Z) = \left\{ 
\begin{array}{rl}    
{\mathcal O}_Z & \quad (i = 0)\\ 
N_{Z/{\mathcal X}} & \quad 
(i = 1) \\
{\mathcal O}(-2)\vert_Z & \quad 
(i = 2) \\ 
0 & \quad \mbox{(otherwise)} 
\end{array}\right.$$ 

We shall prove that 
$H^i(Z, {\underline{\mathrm 
Ext}}^j({\mathcal O}_Z, {\mathcal O}_Z)) 
= 0$ except for $(i,j) = (0,0), 
(1,2)$. For $(i,0)$ with $i > 0$, 
the cohomology clearly vanishes 
because $h^i({\mathcal O}_Z) = 0$ 
for $i > 0$. For $(0,2)$ and 
for $(i,2)$ with $i > 1$, one can 
check that the cohomologies also 
vanish. Hence we only have to prove 
that $h^i(Z, N_{Z/{\mathcal X}}) = 0$ 
for all $i$. By the exact sequence 
above, we immediately see that 
$h^i(Z, N_{Z/{\mathcal X}}) = 0$ for 
$i > 1$. Now let us consider the 
commutative diagram in (3.2.4): 
  
$$ {\tilde Z} \stackrel{\phi}\to Z $$ 
$$ {\tilde p} \downarrow 
\hspace{0.4cm} p \downarrow $$ 
$$ {\tilde N} \stackrel{\bar{\phi}}\to 
\Sigma. $$ 

For $q \in \Sigma$, let $Z_q$ be 
the fiber of $p$ over $q$. If 
$q \ne 0$, then  
$Z_q \cong {\mathbf P}^1$ and 
$N_{Z/{\mathcal X}}\vert_{Z_q} \cong 
{\mathcal O}_{{\mathbf P}^1}(-1)\oplus 
{\mathcal O}_{{\mathbf P}^1}(-1)$. 
Therefore, $H^0(Z, N_{Z/{\mathcal X}}) 
= 0$. Let $p_*{\mathcal O}_Z \to 
R^1p_*({\mathcal O}(-2)\vert_Z)$ be 
the connecting homomorphism 
induced from the sequence above.  
The right hand side can be written as 
    
$R^1p_*({\mathcal O}(-2)\vert_Z) = 
R^1p_*(\phi_*\phi^*({\mathcal O}(-2)\vert_Z)) =$ 

$R^1(p\circ\phi)_*(\phi^*({\mathcal O}(-2)\vert_Z)) 
= R^1({\bar{\phi}}\circ{\tilde p})_*
(\phi^*({\mathcal O}(-2)\vert_Z)) \subset$ 

${\bar{\phi}}_*R^1{\tilde p}_*
(\phi^*({\mathcal O}(-2)\vert_Z)).$ 

Here, in the second equality, we used 
the fact $Z$ has only rational singularities. 
The last map is an inclusion, because 
${\tilde p}_*(\phi^*({\mathcal O}(-2)\vert_Z)) 
= 0$ and hence $R^1{\bar{\phi}}_*
(\phi^*({\mathcal O}(-2)\vert_Z)) = 0$.  
In particular, we see that    
$R^1p_*({\mathcal O}(-2)\vert_Z)$ is a 
torsion free sheaf of rank 1. 
Note that $R^1p_*N_{Z/{\mathcal X}}$ 
is zero outside $0 \in \Sigma$.   
Hence, the connecting homomorphism 
is an isomorphism outside $0$.  
Since $p_*{\mathcal O}_Z$ is reflexive 
and $\dim \Sigma = 6$, we conclude 
that the connecting homomorphism 
is an isomorphism.  Therefore, 
$R^1p_*N_{Z/{\mathcal X}} = 0$ and 
$H^1(Z, N_{Z/{\mathcal X}}) = 0.$   
The results of the lemma follow 
from the local to global spectral 
sequence of Ext.   
\vspace{0.12cm}

{\bf Lemma (4.6.5)}. 

$$\mathrm{Ext}^i({\mathcal O}_{Z'}, 
{\mathcal O}_{{\mathcal X}}) = \left\{ 
\begin{array}{rl}    
0 & \quad (i \ne 3,9)\\ 
H^1(X, {\mathcal O}(-2)\vert_X) & \quad 
(i = 3) \\
{\mathbf C} & \quad 
(i = 9)
\end{array}\right.$$ 

{\em Proof}. Apply 
$\mathrm{Ext}(\bullet, 
{\mathcal O}_{\mathcal X})$ to the 
exact sequence (4.5.1) 
$$ 0 \to {\mathcal O}_M 
\to {\mathcal O}_{Z'} \to 
{\mathcal O}_Z \to 0.$$
We have an exact sequence 
$$  \to \mathrm{Ext}^i
({\mathcal O}_Z, {\mathcal O}_{{\mathcal X}}) 
\to \mathrm{Ext}^i({\mathcal O}_{Z'}, 
{\mathcal O}_{{\mathcal X}}) \to 
\mathrm{Ext}^i({\mathcal O}_M, 
{\mathcal O}_{{\mathcal X}}) \to .$$

 Since $M$ is compact and 
$\omega_{{\mathcal X}}\vert_M 
\cong {\mathcal O}_M$, 
$\mathrm{Ext}^i({\mathcal O}_M, 
{\mathcal O}_{{\mathcal X}})$ is dual 
to $\mathrm{Ext}^{9-i}
({\mathcal O}_{{\mathcal X}}, {\mathcal O}_M) 
= H^{9-i}(M, {\mathcal O}_M).$ 
Hence, we have 
$$\mathrm{Ext}^i({\mathcal O}_M, 
{\mathcal O}_{{\mathcal X}}) = \left\{ 
\begin{array}{rl}    
0 & \quad (i \ne 9)\\ 
{\mathbf C} & \quad 
(i = 9) 
\end{array}\right.$$ 
The exact sequence above 
and (4.6.3) now give the 
result. 
\vspace{0.12cm}

{\bf Lemma (4.6.6)}. 
$$\mathrm{Ext}^i({\mathcal O}_M, 
{\mathcal O}_X) = \left\{ 
\begin{array}{rl}    
0 & \quad (i \ne 8,9)\\ 
{\mathbf C} & \quad 
(i = 8,9)
\end{array}\right.$$ 

$$\mathrm{Ext}^i({\mathcal O}_M, 
{\mathcal O}(2)\vert_X) = 0 \quad 
(\forall i)$$ 
{\em and} 

$$\mathrm{Ext}^i({\mathcal O}_M, 
{\mathcal O}_Z) = \left\{ 
\begin{array}{rl}    
0 & \quad (i \ne 8,9)\\ 
{\mathbf C} & \quad 
(i = 8,9) 
\end{array}\right.$$ 

{\em Proof}. By the exact sequence 
(cf.(4.2)) 
$$ 0 \to {\mathcal O}_{{\mathcal X}} \to 
{\mathcal O}_{{\mathcal X}} \to {\mathcal O}_X 
\to 0$$ 
we have an exact sequence 
$$\mathrm{Ext}^i({\mathcal O}_M, 
{\mathcal O}_{{\mathcal X}}) \stackrel{t}\to 
\mathrm{Ext}^i({\mathcal O}_M, 
{\mathcal O}_{{\mathcal X}}) 
\to 
\mathrm{Ext}^i({\mathcal O}_M, 
{\mathcal O}_X) \to 
\mathrm{Ext}^{i+1}({\mathcal O}_M, 
{\mathcal O}_{{\mathcal X}}),$$
where $t$ is the local coordinate 
of ${\mathbf C}^1$ with $t(0) = 0$ 
(cf.(4.2)). Since $t{\mathcal O}_M = 0$, 
the first map in the sequence is zero. 
On the other hand, since $M$ is 
compact and $\omega_{{\mathcal X}}\vert_M 
\cong {\mathcal O}_M$,  
$\mathrm{Ext}^i({\mathcal O}_M, 
{\mathcal O}_{{\mathcal X}})$ is the dual space of 
$H^{9-i}(M, {\mathcal O}_M)$. Now the first 
claim follows from the exact 
sequence above.    
Next apply $\mathrm{Ext}({\mathcal O}_M, 
\bullet )$ to the exact sequence 
$$ 0 \to {\mathcal O}(2) \to 
{\mathcal O}(2) \to {\mathcal O}(2)\vert_X 
\to 0.$$ Note that 
$\mathrm{Ext}^i({\mathcal O}_M, 
{\mathcal O}(2))$ is the dual space of 
$H^{9-i}(M, {\mathcal O}_M(-2))$ by 
the Serre duality. Since 
$H^i(M, {\mathcal O}_M(-2)) = 0$ for 
all $i$, the second claim follows. 
Finally apply 
$\mathrm{Ext}({\mathcal O}_M, 
\bullet)$ to the exact sequence 
in (4.2): 
$$ 0 \to {\mathcal O}(2)\vert_X 
\to {\mathcal O}_X \to {\mathcal O}_Z 
\to 0. $$ 
By the first and the second claims, 
we have the third statements.  
\vspace{0.12cm}

{\bf Lemma (4.6.7)}.

$$\mathrm{Ext}^i({\mathcal O}_{Z'}, 
{\mathcal O}_Z) = \left\{ 
\begin{array}{rl}    
0 & \quad (i \ne 0,3,8,9)\\ 
H^0(Z, {\mathcal O}_Z) & \quad 
(i = 0) \\
H^1(Z, {\mathcal O}(-2)\vert_Z) & 
\quad (i = 3) \\ 
{\mathbf C} & \quad 
(i = 8,9)
\end{array}\right.$$ 

{\em Proof}. 
By the exact sequence (4.5.1) 
$$ 0 \to {\mathcal O}_M \to 
{\mathcal O}_{Z'} \to 
{\mathcal O}_Z \to 0$$ 
we have an exact sequence 
$$ \to \mathrm{Ext}^i({\mathcal O}_Z, 
{\mathcal O}_Z) \to 
\mathrm{Ext}^i({\mathcal O}_{Z'}, 
{\mathcal O}_Z) \to 
\mathrm{Ext}^i({\mathcal O}_M, 
{\mathcal O}_Z) \to . $$ 
Here we use (4.6.4) and (4.6.6). 
\vspace{0.12cm}

{\bf Lemma (4.6.8)}. 

$$\mathrm{Ext}^5({\mathcal O}_M,  
{\mathcal O}_M) \ne 0.$$
\vspace{0.12cm}

{\em Proof}. Since $M \subset 
{\mathcal X}$ is of locally complete 
intersection, ${\underline{\mathrm{Ext}}}^i
({\mathcal O}_M, {\mathcal O}_M) \cong 
\wedge^i N_{M/{\mathcal X}}$. 
Since $N_{M/X} \cong \Omega^1_M$ 
and $N_{X/{\mathcal X}}\vert_M \cong 
{\mathcal O}_M$, we have an exact 
sequence 
$$ 0 \to \Omega^1_M \to N_{M/{\mathcal X}} \to 
{\mathcal O}_M \to 0. $$ 
This sequence, in particular, yields 
the exact sequence 
$$ 0 \to \Omega^3_M \to 
\wedge^3 N_{M/{\mathcal X}} \to 
\Omega^2_M \to 0.$$ 
The following sequence is exact. 
$$ H^2(\wedge^3 N_{M/{\mathcal X}}) 
\to H^2(M, \Omega^2_M) \to 
H^3(M, \Omega^3_M). $$ 
Since $M \cong G(2,4)$, 
$h^2(M, \Omega^2_M) = 2$ and 
$h^3(M, \Omega^3_M) = 1$. 
Therefore, we have 
$H^2(M, \wedge^3 N_{M/{\mathcal X}}) 
\ne 0.$ 
We use the spectral sequence 
$$ E^{i,j}_2 := H^i({\mathcal X}, 
{\underline{\mathrm{Ext}}}^j({\mathcal O}_M,  
{\mathcal O}_M)) \Rightarrow \mathrm{Ext}^{i+j}
({\mathcal O}_M, {\mathcal O}_M) $$ 
to compute $\mathrm{Ext}^5({\mathcal O}_M, 
{\mathcal O}_M))$. By the argument above, 
$E^{2,3}_2 \ne 0$. Moreover, we can 
check that $E^{2,3}_{\infty} = 
E^{2,3}_2.$ In particular, 
$\mathrm{Ext}^5({\mathcal O}_M, 
{\mathcal O}_M) \ne 0.$   
\vspace{0.12cm}

{\bf Lemma (4.6.9)}. 
$$\mathrm{Ext}^5({\mathcal O}_{Z'}, 
{\mathcal O}_M) \ne 0,$$ {\em and}
$$\mathrm{Ext}^5({\mathcal O}_{Z'}, 
{\mathcal O}_{Z'}) \ne 0.$$ 
{\em Proof}. By the exact sequence 
(4.5.1) 
$$ 0 \to {\mathcal O}_M \to 
{\mathcal O}_{Z'} \to {\mathcal O}_Z 
\to 0$$
we have an exact sequence
$$ \mathrm{Ext}^5({\mathcal O}_Z, 
{\mathcal O}_M) \to 
\mathrm{Ext}^5({\mathcal O}_{Z'}, 
{\mathcal O}_M) \to 
\mathrm{Ext}^5({\mathcal O}_M, 
{\mathcal O}_M) \to 
\mathrm{Ext}^6({\mathcal O}_Z, 
{\mathcal O}_M).$$ 
By the Serre duality, 
$\mathrm{Ext}^i({\mathcal O}_Z, 
{\mathcal O}_M) \cong 
(\mathrm{Ext}^{9-i}({\mathcal O}_M, 
{\mathcal O}_Z))^*$. For $i = 5,6$, 
these vanish by (4.6.6). 
Now the first claim follows from (4.6.8). 
We apply $\mathrm{Ext}({\mathcal O}_{Z'}, 
\bullet)$ to the exact sequence (4.5.1) 
above. Since $\mathrm{Ext}^i({\mathcal O}_{Z'}, 
{\mathcal O}_Z) \ne 0$ for $i = 4,5$ by 
(4.6.7), we have an isomorphism 
$\mathrm{Ext}^5({\mathcal O}_{Z'}, 
{\mathcal O}_M) \cong \mathrm{Ext}^5({\mathcal O}_{Z'}, 
{\mathcal O}_{Z'}).$ The second claim 
follows from the first claim. 
\vspace{0.12cm}

{\bf Lemma (4.6.10)}. $${\mathrm 
Ext}^5(I_{Z'}, I_{Z'}) \ne 0.$$ 
{\em Proof}. In the exact sequence 
$$\mathrm{Ext}^5({\mathcal O}_{{\mathcal X}}, 
I_{Z'})  \to 
\mathrm{Ext}^5(I_{Z'}, 
I_{Z'}) \to 
\mathrm{Ext}^6({\mathcal O}_{Z'}, 
I_{Z'}) \to 
\mathrm{Ext}^6({\mathcal O}_{{\mathcal X}}, 
I_{Z'})$$
the first term and the last term 
both vanish by (4.6.1). Hence 
$$ \mathrm{Ext}^5(I_{Z'}, 
I_{Z'}) \cong 
\mathrm{Ext}^6({\mathcal O}_{Z'}, 
I_{Z'}).$$  
In the exact sequence 
$$\mathrm{Ext}^5({\mathcal O}_{Z'}, 
{\mathcal O}_{{\mathcal X}})  \to 
\mathrm{Ext}^5({\mathcal O}_{Z'}, 
{\mathcal O}_{Z'}) \to 
\mathrm{Ext}^6({\mathcal O}_{Z'}, 
I_{Z'}) \to 
\mathrm{Ext}^6({\mathcal O}_{Z'}, 
{\mathcal O}_{{\mathcal X}})$$ 
the first and the last terms 
both vanish by (4.6.5). 
Hence 
$$ \mathrm{Ext}^5({\mathcal O}_{Z'}, 
{\mathcal O}_{Z'}) \cong 
\mathrm{Ext}^6({\mathcal O}_{Z'}, 
I_{Z'}). $$ 
By (4.6.9), 
$$ \mathrm{Ext}^5({\mathcal O}_{Z'}, 
{\mathcal O}_{Z'}) \ne 0.$$
\vspace{0.12cm}
 
{\bf  Observation (4.7).} 
 {\em $\Phi$ is not fully 
faithful.}  
\vspace{0.12cm}

{\em Proof}. By (4.6.10) $\mathrm{Ext}^5
(\Phi({\mathcal O}^+(1), \Phi({\mathcal O}^+(1))) 
\ne 0$. 
Since $\mathrm{Ext}^i
({\mathcal O}^+(1), {\mathcal O}^+(1)) 
= 0$ for all $i > 0$, this implies 
that $\Phi$ is not fully faithful.  
\vspace{0.2cm}

{\bf Lemma (4.8)}. $\Phi = \Psi$ 
{\em and } $\Phi_0 = \Psi_0.$  
\vspace{0.12cm}

{\em Proof}.  
We put ${\hat{\mathcal X}} := 
E(H) \times_{{\bar{\mathcal N}}}E(H^*)$ 
and ${\hat X} := 
T^*G \times_{{\bar N}}T^*G^+.$ 
Let $\alpha: 
{\tilde{\mathcal X}} \to {\hat{\mathcal X}}$ 
and $\alpha_0: {\tilde X} \to 
{\hat X}$ 
be the natural morphisms 
(cf.(3.7)). 
Note that $\nu = \mu \circ \alpha$ 
and $\nu_0 = \mu_0 \circ \alpha_0.$ 
Since ${\hat{\mathcal X}}$ has only 
rational singularities by (3.7), we can write 
$$ {\mathbf R}\nu_* \circ 
{\mathbf L}(\nu^+)^* = 
{\mathbf R}\mu_*\circ {\mathbf R}\alpha_* 
\circ {\mathbf L}\alpha^* \circ 
{\mathbf L}(\mu^+)^* = 
{\mathbf R}\mu_* \circ {\mathbf L}(\mu^+)^*. $$  
Therefore, $\Phi = \Psi.$  
We next claim that 
${\mathbf R}(\alpha_0)_*{\mathcal O}_{{\tilde X}} 
= {\mathcal O}_{{\hat X}}.$  
Since $R^1\alpha_*{\mathcal O}_{{\tilde{\mathcal X}}} 
= 0$, we have the exact sequence 
$$ 0 \to \alpha_*{\mathcal O}_{{\tilde{\mathcal X}}} 
\stackrel{t}\to \alpha_*{\mathcal O}_{{\tilde{\mathcal X}}} \to 
\alpha_*{\mathcal O}_{{\tilde X}} \to 0.$$ 
Since $\alpha_*{\mathcal O}_{{\tilde{\mathcal X}}} = 
{\mathcal O}_{{\hat{\mathcal X}}}$, this sequence is 
identified with 
$$ 0 \to {\mathcal O}_{{\hat{\mathcal X}}} \stackrel{t}\to 
{\mathcal O}_{{\hat{\mathcal X}}} \to {\mathcal O}_{{\hat X}} \to 
0.$$  Hence, $\alpha_*{\mathcal O}_{{\tilde X}} 
\cong {\mathcal O}_{{\hat X}}$. Therefore,  
$(\alpha_0)*{\mathcal O}_{{\tilde X}} = 
{\mathcal O}_{{\hat X}}.$ 
For $i > 0$, the exact sequence 
$$R^i\alpha_*{\mathcal O}_{{\tilde{\mathcal X}}} 
\to R^i(\alpha_0)_*{\mathcal O}_{{\tilde X}} 
\to R^{i+1}\alpha_*{\mathcal O}_{{\tilde{\mathcal X}}} 
$$ yields $R^i(\alpha_0)_*{\mathcal O}_{{\tilde X}} 
= 0$. Our claim is now justified.   
Then we can write 
$$ {\mathbf R}(\nu_0)_* \circ 
{\mathbf L}(\nu^+_0)^* = 
{\mathbf R}(\mu_0)_*\circ {\mathbf R}(\alpha_0)_* 
\circ {\mathbf L}(\alpha_0)^* \circ 
{\mathbf L}(\mu^+_0)^* = 
{\mathbf R}(\mu_0)_* \circ {\mathbf L}(\mu^+_0)^*. $$ 
Therefore, $\Phi_0 = \Psi_0.$   
\vspace{0.2cm}

{\bf Observation (4.9)}. 
(1). {\em $\Psi$ is not fully faithful.} 

(2) {\em $\Psi_0$ is not an equivalence.}
\vspace{0.12cm}

{\em Proof}. (1): This is clear from (4.7) and 
(4.8). 

(2): (i) The functor 
${\Phi}' := {\mathbf R}
{{\nu}^+}_*({\mathbf L}\nu^*(\bullet) \otimes 
\omega_{{\tilde{\mathcal X}}/{\mathcal X}})$  
is the right adjoint of $\Phi$, where 
$\omega_{{\tilde{\mathcal X}}/{\mathcal X}} 
:= \omega_{{\tilde{\mathcal X}}} \otimes 
\nu^*(\omega_{{\mathcal X}}^{-1})$.  
On the other hand, the functor 
${\Phi}'' := {\mathbf R}
{{\nu}^+}_*({\mathbf L}\nu^*(\bullet) \otimes 
\omega_{{\tilde{\mathcal X}}/{\mathcal X}^+})$ 
is the left adjoint of $\Phi$.  
But, since $\nu^*\omega_{{\mathcal X}} 
= (\nu^+)^*\omega_{{\mathcal X}^+}$, 
these functors coincides: 
${\Phi}' = {\Phi}''.$ 
Now, ${\Phi}'_0 := 
{\mathbf R}
{{\nu}_0^+}_*({\mathbf L}(\nu_0)^*(\bullet) \otimes 
\omega_{{\tilde X}/X})$
becomes the adjoint of $\Phi_0$. 

(ii) By the projection formula, 
$\Phi'$ coincides with the functor 
$$\Psi' := {\mathbf R}(\mu^+)_*({\mathbf L}
\mu^*(\bullet)\otimes {\mathbf R}\alpha_*
\omega_{{\tilde{\mathcal X}}/{\mathcal X}}).$$ 
Similarly, ${\Phi}'_0$ coincides with 
the functor
$${\Psi}'_0 :=  
{\mathbf R}((\mu_0)^+)_*({\mathbf L}
(\mu_0)^*(\bullet)\otimes {\mathbf R}(\alpha_0)_*
\omega_{{\tilde X}/X}).$$
By the Grauert-Riemmenschneider vanishing, 
we have ${\mathbf R}\alpha_*
\omega_{{\tilde{\mathcal X}}/{\mathcal X}} 
= \alpha_*\omega_{{\tilde{\mathcal X}}/{\mathcal X}}$. 
By the exact sequence 
$$ R^i\alpha_*\omega_{{\tilde{\mathcal X}}/
{\mathcal X}} \stackrel{t}\to 
R^i\alpha_*\omega_{{\tilde{\mathcal X}}/
{\mathcal X}} \to R^i(\alpha_0)_*
\omega_{{\tilde X}/X} \to 
R^{i+1}\alpha_*\omega_{{\tilde{\mathcal X}}/
{\mathcal X}}$$ 
we see that 
${\mathbf R}(\alpha_0)_*
\omega_{{\tilde X}/X} 
= (\alpha_0)_*\omega_{{\tilde X}/X}$ and 
$(\alpha_0)_*\omega_{{\tilde X}/X} \cong 
\alpha_*\omega_{{\tilde{\mathcal X}}/{\mathcal X}}
\otimes_{{\mathcal O}_{{\hat{\mathcal X}}}}
{\mathcal O}_{{\hat X}}$. Let $j: {\hat X} \to 
{\hat{\mathcal X}}$ be the inclusion map. 
Since $\alpha_*\omega_{{\tilde{\mathcal X}}/{\mathcal X}}$ 
is flat over ${\mathbf C}^1$ by (3.7), 
${\mathbf L}j^*\alpha_*\omega_{{\tilde{\mathcal X}}/{\mathcal X}} 
= (\alpha_0)_*\omega_{{\tilde X}/X}$.  
Therefore, 
$$(4.9.1) \qquad \Psi' = {\mathbf R}(\mu^+)_*({\mathbf L}
\mu^*(\bullet)\otimes \alpha_*
\omega_{{\tilde{\mathcal X}}/{\mathcal X}}),$$ 
$$(4.9.2) \qquad {\Psi}'_0 =  
{\mathbf R}((\mu_0)^+)_*({\mathbf L}
(\mu_0)^*(\bullet)\otimes 
{\mathbf L}j^*\alpha_*\omega_{{\tilde{\mathcal X}}/{\mathcal X}}).$$
Since $\Psi = \Phi$ and 
$\Psi_0 = \Phi_0$ by (4.8), $\Psi'$ and 
${\Psi}'_0$ are adjoint of $\Psi$ and 
$\Psi_0$ by (i). Let $i: X^+ 
\to {\mathcal X}^+$ be the inclusion map. 
By (4.9.1) and (4.9.2) 
we can apply [Ch, Lemma 6.1, Ka 2, Lemma 5.6] 
to conclude that  
the following diagram commutes  
$$ D(X^+) \stackrel{{\Psi}'_0 \circ {\Psi}_0}
\to D(X^+) $$ 
$$ {\mathbf R}i_* \downarrow 
\hspace{0.8cm}
{\mathbf R}i_* \downarrow $$ 
$$ D({\mathcal X}^+)   
\stackrel{{\Psi}' \circ \Psi}\to 
D({\mathcal X}^+).$$

(iii) Assume now that $\Psi_0$ is 
an equivalence. Then the quasi-inverse 
of $\Psi_0$ coincides with the adjoint 
${\Psi}'_0$. In particular, 
${\Psi}'_0\circ {\Psi}_0 \cong 
id_{D(X^+)}$.  
Let $\Omega := 
\{{\mathcal O}_p\}_{p \in {\mathcal X}^+},$ 
where ${\mathcal O}_p$ are structure 
sheaves of the closed points of 
${\mathcal X}^+$.  Then $\Omega$ 
is the spanning class for 
$D({\mathcal X}^+)$(cf.
[Br 2], Example 2.2). It is clear that 
${\Psi}' \circ \Psi ({\mathcal O}_p) = 
{\mathcal O}_p$ for $p \in 
{\mathcal X}^+\setminus X^+.$  
For $p \in X^+$, since 
${\mathbf R}i_*{\mathcal O}_p = 
{\mathcal O}_p$, we have 
${\Psi}' \circ \Psi ({\mathcal O}_p) = 
{\mathcal O}_p$ by the commutative 
diagram because 
${\Psi}'_0 \circ {\Psi}_0  
\cong id_{D(X^+)}$.  
Then, by the next lemma, we 
see that $\Psi$ is fully faithful; 
but this contradicts (1) .   
\vspace{0.12cm}

{\bf Lemma (4.10)} ([Ka 1, Lemma (5.4)]). 
{\em Let $f: A \to B$ be a functor with 
the right adjoint $g$ and the left adjoint 
$h$. Let $\Omega$ be a spanning class 
for $A$. Assume that $\forall \omega 
\in \Omega$, $\omega \cong 
g \circ f (\omega)$ and 
$h \circ f (\omega) \cong \omega.$ 
Then, $g\circ f (a) \cong h\circ f(a) 
\cong a$ for all $a \in A$. Moreover, 
$f$ is fully faithful.}  
\vspace{0.12cm}

{\em Proof}. For $a \in A$ and $\omega 
\in \Omega$, we have:  
\vspace{0.12cm}

$ \mathrm{Hom}(h \circ f (a), \omega) 
\cong \mathrm{Hom}(f(a), f(\omega))  
\cong $

$\mathrm{Hom}(a, g \circ f (\omega)) 
\cong \mathrm{Hom}(a, \omega), $  
\vspace{0.12cm}

$\mathrm{Hom}(\omega, g \circ f (a)) 
\cong \mathrm{Hom}(f(\omega), f(a)) 
\cong$ 

$\mathrm{Hom}(h \circ f (\omega), a) 
\cong \mathrm{Hom}(\omega, a).$  
\vspace{0.12cm}

Therefore, $a \cong g \circ f(a)$ and 
$h \circ f (a) \cong a.$
Since 
$$ \mathrm{Hom}(f(a), f(b)) \cong 
\mathrm{Hom}(h \circ f (a), b) 
\cong \mathrm{Hom}(a, b), $$ 
$f$ is fully faithful.
\vspace{0.2cm}

\begin{center}
{\bf \S 5. Complete flag varieties}
\end{center}

(5.1). Let $H$ be a ${\mathbf C}$-vector 
space of $\dim = h$. We denote by 
$F$ the complete flag 
vatiety. Namely, 
$$ F := \{V_1 \subset V_2 \subset ... \subset 
V_{h-1} \subset H ; \dim V_i = i \}.$$       
Let 
$$\tau_1 \subset \tau_2 \subset ... 
\subset \tau_{h-1} \subset H \otimes_{{\mathbf C}}
{\mathcal O}_F$$ 
be the universal sub-bundles. We put $q_i 
:= H \otimes_{{\mathbf C}}{\mathcal O}_F/\tau_i$ 
and call them the universal quotient bundles. 
\vspace{0.12cm}

(5.2). Let $\pi: T^*F \to F$ be the cotangent bundle 
of $F$. The nilpotent variety ${\bar N}(H)$ 
is defined as 
$$ {\bar N}(H) := \{A \in \mathrm{End}(H); 
A^h = 0 \}. $$ 
A point of $T^*F$ is expressed as a pair 
$(p, \phi)$ of $p \in F$ and $\phi \in 
\mathrm{End}(H)$ such that 
$$ \phi(H) \subset \tau_{h-1}(p), 
\phi(\tau_{h-1}(p)) \subset 
\tau_{h-2}(p), ..., \phi(\tau_1(p)) = 0.$$ 
The Springer resolution 
$$ s: T^*F \to {\bar N}(H)$$ 
is defined as $s((p, \phi)) := \phi$.    
The Springer resolution $s$ has the 
following properties. 

(5.2.1). $s^{-1}(0) = F$. 
\vspace{0.12cm}

(5.2.2). Let  
${\bar N}(H)_{sing}$ be the 
singular locus of ${\bar N}(H)$. 
Then 
$${\bar N}(H)_{sing} := 
\{A \in {\bar N}(H); \mathrm{rank}
(A) \leq h-2\}. $$ 
Let ${\bar N}(H)^0_{sing}$ be the 
open orbit consisting of the matrices 
conjugate to 
$$ \left(
\begin{array}{cccccc}
0 & 0 & .. & .. & .. & ..\\
0 & 0 & 1 & 0 & .. & .. \\ 
0 & 0 & 0 & 1 & 0 & .. \\
.. & .. & .. & .. & 1 & .. \\ 
.. & .. & .. & .. & .. & 1 \\ 
0 & 0 & .. & .. & .. & 0  
\end{array}\right). $$ 
For $A \in {\bar N}(H)^0_{sing}$, 
$s^{-1}(A)$ is a tree of ${\mathbf P}^1$ 
with the $A_{h-1}$-configuration.    
\vspace{0.12cm}

(5.2.3). For $A \in {\bar N}(H)\setminus 
{\bar N}(H)_{sing}$, $s^{-1}(A)$ is one 
point.  
\vspace{0.12cm}

(5.3). The complete flag variety $F$ 
has $h-1$ natural fibrations 
$f_1: F \to F(2,3, ..., h-1,H)$, 
$f_2: F \to 
F(1,3, ..., h-1, H)$ ..., and 
$f_{h-1}: F \to F(1,2, ..., h-2, H)$. 
Each fibration is a ${\mathbf P}^1$ 
bundle. Assume that $f_i = \Phi_{\vert L_i \vert}$ 
for a line bundle $L_i$ on $F$. Then 
$\pi^*L_i \in \mathrm{Pic}(T^*F)$ defines 
a birational morphism 
$$s_i: T^*F \to X_i$$ 
over ${\bar N}(H)$, where $X_i$ is a normal 
variety which factorize $s$ as $T^*F 
\stackrel{s_i}\to X_i \to {\bar N}(H)$.  
Over $0 \in {\bar N}(H)$, $s_i$ restricts 
to the fibration $f_i : F \to 
F(1,2, ..., i-1, i+1, ... h-1,H)$. 
Let $E_i$ be the exceptional locus of 
$s_i$. Then $E_i \to s_i(E_i)$ is a 
${\mathbf P}^1$ bundle. In $T^*F$, each fiber 
of this ${\mathbf P}^1$-bundle has the 
normal bundle ${\mathcal O}_{{\mathbf P}^1}
^{\oplus 2\dim F - 2}
\oplus {\mathcal O}_{{\mathbf P}^1}(-2).$  
Note that, in the family of 
rational curves: $E_i \to s_i(E_i)$, 
each fiber of $f_i$ deforms to 
one of ${\mathbf P}^1$ in the tree of (5.2.2).  
\vspace{0.12cm}

(5.4). We put $sl(H) := \{A \in \mathrm{End}(H); 
tr(A) = 0\}.$ For $A \in \mathrm{End}(H)$, 
let $\phi_A(x)$ be the characteristic polynomial 
of $A$: 
$$ \phi_A(x) := \mathrm{det}(xI - A). $$ 
Let $\phi_i(A)$ be the coefficient of $x^{h-i}$ 
in $\phi_A(x)$. If $A \in sl(H)$, then $\phi_1(A) 
= 0$. We define the {\em characteristic map}  
$$ {\mathrm ch}: sl(H) \to {\mathbf C}^{h-1} $$ 
as ${\mathrm ch}(A) := (\phi_2(A), ..., \phi_h(A)).$ 
Note that ${\mathrm ch}^{-1}(0) = {\bar N}(H).$ 
\vspace{0.12cm}

(5.5). We shall define a simultaneous 
resolution of ${\mathrm ch}: sl(H) \to 
{\mathbf C}^{h-1}$ up to a finite cover. 
First, we shall define a vector bundle 
${\mathcal E}(H)$ over $F$ and an exact sequence 
$$ 0 \to T^*F \to {\mathcal E}(H) \stackrel{\eta}\to 
{\mathcal O}_F^{\oplus h-1} \to 0. $$ 
Let $T^*F(p)$ be the cotangent space of 
$F$ at $p \in F$. Then, by (5.2), for a suitable 
basis of $H$, $T^*F(p)$ consists of the matrices 
of the following form 
$$ \left(
\begin{array}{ccccc}
0 & * & .. & .. & .. \\
0 & 0 & * & .. & .. \\ 
.. & .. & .. & .. & .. \\
.. & .. & .. & 0 & * \\ 
0 & 0 & .. & 0 & 0   
\end{array}\right). $$ 
Let ${\mathcal E}(H)(p)$ be the vector 
subspace of $sl(H)$ consisting of 
the matrices $A$ of the following 
form 
$$ \left(
\begin{array}{ccccc}
\alpha_1 & * & .. & .. & .. \\
0 & \alpha_2 & * & .. & .. \\ 
.. & .. & .. & .. & .. \\
.. & .. & .. & \alpha_{h-1} & * \\ 
0 & 0 & .. & 0 & \alpha_h   
\end{array}\right).$$   
Here $\alpha_1 + ... + \alpha_h 
= 0.$ 
We define a map $\eta(p): {\mathcal E}(H)(p) 
\to {\mathbf C}^{\oplus h-1}$ as 
$\eta(p)(A) := (\alpha_1, ..., \alpha_{h-1}).$ 
Then we have an exact sequence of 
vector spaces 
$$ 0 \to T^*F(p) \to {\mathcal E}(H)(p) 
\stackrel{\eta(p)}\to {\mathbf C}^{\oplus 
h-1} \to 0. $$  
We put ${\mathcal E}(H) := \cup_{p \in F}
{\mathcal E}(H)(p)$. Then ${\mathcal E}(H)$ 
becomes a vector bundle over $F$, and 
we get the desired exact sequence. 
Each point of ${\mathcal E}(H)$ is 
expressed as a pair of $p \in F$ 
and $\phi \in {\mathcal E}(H)(p)$. Now 
we define 
$$ {\tilde s}: {\mathcal E}(H) \to 
sl(H)$$ 
as ${\tilde s}((p, \phi)) := \phi.$ 
The birational morphisms 
$s_i: T^*G \to X_i$ extend to the 
birational morphisms 
$${\tilde s}_i: 
{\mathcal E}(H) \to {\mathcal X}_i$$
and 
${\tilde s}$ is factorized as 
${\mathcal E}(H) 
\to {\mathcal X}_i \to sl(H).$    
Let ${\mathcal E}_i$ be the exceptional locus of 
${\tilde s}_i$. Then ${\mathcal E}_i \to 
{\tilde s}_i({\mathcal E}_i)$ is a 
${\mathbf P}^1$ bundle. In ${\mathcal E}(H)$, each fiber 
of this ${\mathbf P}^1$-bundle has the 
normal bundle ${\mathcal O}_{{\mathbf P}^1}
^{\oplus 2\dim F + h - 4}
\oplus {\mathcal O}_{{\mathbf P}^1}(-1)
\oplus {\mathcal O}_{{\mathbf P}^1}(-1).$ 
Let ${\tilde \eta}: {\mathcal E}(H) \to 
{\mathbf C}^{h-1}$ be the morphism 
induced by $\eta$. We define a finite 
Galois cover 
$$\varphi: {\mathbf C}^{h-1} 
\to {\mathbf C}^{h-1}$$ 
as 
$\varphi (\alpha_1, ..., \alpha_{h-1}) 
:= (\phi_2(A), ..., \phi_h(A))$  
where $A$ is the diagonal matrix 
$$\left(
\begin{array}{ccccc}
\alpha_1 & 0 & .. & .. & .. \\
0 & \alpha_2 & 0 & .. & .. \\ 
.. & .. & .. & .. & .. \\
.. & .. & .. & \alpha_{h-1} & 0 \\ 
0 & 0 & .. & 0 & \alpha_h   
\end{array}\right)$$ 
with $\alpha_h = -\alpha_1 - ... - 
\alpha_{h-1}$.  
Then we have a commutative diagram 

$$ {\mathcal E}(H) \stackrel{{\tilde s}}\to 
sl(H)$$ 
$$ {\tilde \eta} \downarrow 
\hspace{0.4cm} 
{\mathrm ch} \downarrow $$
$$ {\mathbf C}^{h-1} \stackrel{\varphi}\to 
{\mathbf C}^{h-1}.$$ 
Let ${\mathrm ch}': 
sl(H)\times_{{\mathbf C}^{h-1}}
{\mathbf C}^{h-1} \to {\mathbf C}^{h-1}$ 
be the pull-back of ${\mathrm ch}$ 
by $\varphi$. By the commutative 
diagram, we have a morphism 
$$\beta: {\mathcal E}(H) \to  
sl(H)\times_{{\mathbf C}^{h-1}}
{\mathbf C}^{h-1}.$$ 
By this morphism, ${\tilde{\eta}}: 
{\mathcal E}(H) \to {\mathbf C}^{h-1}$ becomes 
a simultaneous resolution of ${\mathrm 
ch}'$. 
\vspace{0.12cm}

(5.6). For the dual space 
$H^*$, we define the complete 
flag variety $F^+$. We denote 
by $T^*F^+$ the cotangent 
bundle of $F^+$. We define a 
nilpotent variety ${\bar N}(H^*)$ 
and the Springer resolution 
$s^+: T^*F^+ \to {\bar N}(H^*)$ 
in the same way as (5.2). 
Let ${\tilde{\eta}}^+: {\mathcal E}(H^*) 
\to {\mathbf C}^{h-1}$ be the corresponding 
objects of (5.5) for $H^*$. 
The natural isomorphism 
$\mathrm{End}(H) \cong 
\mathrm{End}(H^*)$ induces 
an isomorphism 
$\iota: sl(H) \cong sl(H^*).$ 
Let ${\mathrm ch}^+: sl(H^*) 
\to {\mathbf C}^{h-1}$ be the 
characteristic map. Then 
$\iota$ is compatible with 
${\mathrm ch}$ and ${\mathrm ch}^+$: 
${\mathrm ch}^+ \circ \iota = 
{\mathrm ch}.$  We have an 
isomorphism 
$$ sl(H)\times_{{\mathbf C}^{h-1}}
{\mathbf C}^{h-1} \stackrel{\iota 
\times id}\to 
sl(H^*)\times_{{\mathbf C}^{h-1}}
{\mathbf C}^{h-1}.$$ 
This induces a birational map 
$$f: {\mathcal E}(H) - - \to  
{\mathcal E}(H^*).$$ 

(5.7). Let $L^+ \in \mathrm{Pic}
({\mathcal E}(H))$ be a $\beta^+$-ample 
line bundle, where 
$\beta^+: {\mathcal E}(H^*) 
\to sl(H^*)\times_{{\mathbf C}^{h-1}}
{\mathbf C}^{h-1}$ is the simultaneous 
resolution defined in (5.5). Denote by 
$L \in \mathrm{Pic}({\mathcal E}(H))$ its 
proper transform by $f$. For $\sigma 
\in {\mathrm Gal}(\varphi)\cong {\mathrm 
S}_h$, we consider the isomorphism 
$sl(H)\times_{{\mathbf C}^{h-1}}
{\mathbf C}^{h-1} \stackrel{id \times 
\sigma}\to sl(H)\times_{{\mathbf C}^{h-1}}
{\mathbf C}^{h-1}$. This isomorphism induces 
a birational map 
$\phi_{\sigma}: {\mathcal E}(H) - - \to 
{\mathcal E}(H)$. \vspace{0.12cm}

{\bf Lemma (5.7.1)}. {\em For a suitable 
$\sigma \in {\mathrm Gal}(\varphi)$, the 
proper transform of $L$ by $\phi_{\sigma}$ 
becomes $\beta$-ample. In particular, for 
this $\sigma \in {\mathrm Gal}(\varphi)$,  
the composite $f_{\sigma} := f\circ \phi_{\sigma}: 
{\mathcal E}(H) - - \to {\mathcal E}(H^*)$ becomes an 
isomorphism.}       
\vspace{0.12cm}

{\em Proof}. Note that ${{\mathrm ch}'}^{-1}
(0) = {\bar N}(H)$. For $p \in 
{\bar N}(H)^0_{sing}$ (cf.(5.2.2)), 
$\beta^{-1}(p)$ is a tree of ${\mathbf P}^1$ 
with $A_{h-1}$-configuration. Let $C_1, ..., 
C_{h-1}$ be the irreducible components of 
this tree. By [Re, \S 7], we can take $\sigma 
\in {\mathrm Gal}(\varphi)$ in such a way 
that the proper 
transform $L_{\sigma}$ of $L$ by $\phi_{\sigma}$ 
has positive intersections with all 
$C_i$. We next observe the central fiber 
$\beta^{-1}(0) = F$. 
As in (5.3), $F$ has $h-1$ different 
${\mathbf P}^1$ fibrations $f_1$, ..., $f_{h-1}$. 
Let $l_j$ be a fiber of $f_j$.  
The cone $NE(F)$ of 
effective 1-cycles is a polyhedral cone 
generated by $h-1$ rays ${\mathbf R}^+[l_j]$. 
Since $l_i$ deforms to one of $C_i$'s in 
${\mathcal E}(H)$, we see that $(L_{\sigma}. l_j)
> 0$ for all $j$. This implies that $L_{\sigma}$ 
is $\beta$-ample.  
\vspace{0.12cm}

(5.8). ${\mathrm Gal}(\varphi) 
= \{\phi_{\sigma}\}$ contains 
the Atiyah flops along ${\mathcal E}_i 
\subset {\mathcal E}(H)$ for $i = 1, ..., h-1$. 
They are generators of  
${\mathrm Gal}(\varphi)$. 
We denote them 
by $\phi_i$. 
One can choose an isomorphism 
${\mathrm Gal}(\varphi) \cong S_h$ 
in such a way that 
$\phi_i$ is sent to $(i,i+1)$.  
For the flop $\phi_i$: 
$$ {\mathcal E}(H) \rightarrow {\mathcal X}_i 
\leftarrow {\mathcal E}(H),$$ 
let $$\Psi_i: D({\mathcal E}(H)) \to 
D({\mathcal E}(H))$$ be the functor defined 
by the fiber product 
${\mathcal E}(H)\times_{{\mathcal X}_i}
{\mathcal E}(H)$. 
\vspace{0.12cm}

{\bf Proposition (5.8.1)}. 
{\em $\Psi_i$ is an equivalence.}
\vspace{0.12cm}

{\em Proof}. Let $\Psi'_i : 
D({\mathcal E}(H)) \to 
D({\mathcal E}(H))$ be the adjoint 
functor of $\Psi_i$ (cf. the 
proof of (4.9),(2)).   
Let $\Omega := \{{\mathcal O}_p\}_{p \in 
{\mathcal E}(H)}$. Then $\Omega$ is 
a spanning class for $D({\mathcal E}(H))$.
We show that $\Psi_i \circ \Psi'_i (\omega) 
\cong \omega$ and $\Psi'_i \circ 
\Psi_i (\omega) \cong \omega$ for all 
$\omega \in \Omega$. Then, by Lemma (4.10), 
we conclude that $\Psi_i \circ \Psi'_i \cong 
id$ and $\Psi'_i \circ \Psi_i \cong id$.    
The problem being local, we can 
replace ${\mathcal E}(H)$ by ${\mathcal X} :=  
X \times S$ with $X$ a smooth quasi-projective 
threefold 
containing $(-1,-1)$-curve $C$ and 
with $S$ a smooth quasi-projective 
variety. Let ${\bar X}$ be the threefold 
obtained from $X$ by contracting $C$ to 
a point. Let 
$$ X \rightarrow {\bar X} \leftarrow X^+$$ 
be the Atiyah flop along $C \subset X$, and 
let 
$$ {\mathcal X} \rightarrow {\bar X}\times 
S \leftarrow {\mathcal X}^+$$  
be the product of the Atiyah flop with $S$, 
where ${\mathcal X}^+ 
:= X^+ \times S$. Let $\Psi : D({\mathcal X}) 
\to D({\mathcal X}^+)$ be the functor defined 
by the fiber product 
${\mathcal X}\times_{{\bar X}\times S}
{\mathcal X}^+$. Denote by $\Psi'$ its adjoint.  
Moreover, let $\Psi_0: D(X) 
\to D(X^+)$ be the functor defined by 
$X \times_{{\bar X}}X^+$. We already know 
that $\Psi_0$ is an equivalence (cf. [Na, Ka 2]). 
Then, by the same argument as (4.9),(2), we 
see that $\Psi \circ \Psi'({\mathcal O}_p) 
\cong {\mathcal O}_p$ for $p \in {\mathcal X}^+$ 
and $\Psi' \circ \Psi ({\mathcal O}_p) \cong 
{\mathcal O}_p$ for $p \in {\mathcal X}$.    
\vspace{0.15cm}
 
(5.9). For $\sigma \in {\mathrm Gal}(\varphi)$, 
the birational automorphism $\phi_{\sigma}$ 
can be decomposed into a finite sequence of 
Atiyah flops $\phi_i$. 
By (5.7.1), $f = 
f_{\sigma}\circ \phi_{\sigma^{-1}}$ 
for a suitable $\sigma \in 
{\mathrm Gal}(\varphi).$ 
The $\sigma$ corresponds to 
$(12)(23)...(h-2,h-1)(h-1,h)(h-2,h-1)...(12)$ 
by the identification ${\mathrm Gal}(\varphi) 
\cong S_h$ in (5.8). 
Now $\phi_{\sigma^{-1}}$ is decomposed as  
$$ \phi_{\sigma^{-1}} = \phi_1\circ ... \circ 
\phi_{h-2}\circ\phi_{h-1}\circ\phi_{h-2}\circ ... 
\circ \phi_1. $$ 
\vspace{0.12cm} 

{\bf Theorem (5.9.1)}. {\em 
$f$ induces an equivalence of 
derived categories} 
$$ \Psi_f := \Psi_{f_{\sigma}}\circ 
\Psi_1\circ ... \circ 
\Psi_{h-2}\circ\Psi_{h-1}\circ\Psi_{h-2}\circ ... 
\circ \Psi_1: D({\mathcal E}(H)) 
\to D({\mathcal E}(H^*)).$$ 
{\em Moreover, $\Psi_f$ induces an 
equivalence} 
$$ (\Psi_f)_0: D(T^*F) \to D(T^*F^+).$$       

{\em Proof}. The first claim is clear from 
(5.8.1). The second claim follows from [Ka 2, 
Lemma 5.6, Corollary 5.7].     
\vspace{0.2cm}

{\em Acknowledgement}: The author thanks H. Nakajima 
for a helpful letter concerning the similar question 
to [Na] for stratified Mukai flops. He also thanks 
E. Markman for providing him with the proof of 
(2.7.1), which enables him to remove an assumption 
from Theorem (2.7) of the first version.

\vspace{0.2cm}

\begin{center}
Department of Mathematics, Graduate School of 
Science, Kyoto University, 

Kita-shirakawa Oiwake-cho, Kyoto, 
606-8502, Japan 
\end{center}

\end{document}